\documentclass[a4paper,11pt]{article}
\usepackage[a4paper, margin=1.2in]{geometry}
\usepackage{mathrsfs,fancybox,amsfonts,latexsym,theorem,amssymb,amsmath,newlfont, epsf}
\usepackage{graphicx}
\usepackage{multimedia}
\usepackage{pgf}
\usepackage{fancyhdr}
\usepackage{times}
\usepackage{pstricks,pstricks-add,pst-math,pst-xkey} 
\usepackage{nicefrac}
\usepackage{bbold}

\usepackage{srcltx}
\theoremstyle{plain}

\newtheorem{thm}{Theorem}[section]

\newtheorem{lem}[thm]{Lemma}

\def\PP{{\mathbb P}}
\def\RR{{\mathbb R}}

\def\FF{{\cal F}}

\def\<{\langle}
\def\>{\rangle}

\title{Recovering the uniform boundary observability with spectral Legendre-Galerkin formulations of the 1-D wave equation}\author{\Large Ludovick Gagnon$^{1}$ and Jos\'e M. Urquiza$^{2}$\\$^{1}$ Laboratoire J.A. Dieudonn\'e, Universit\'e C\^ote d'Azur,\\ Nice, France. gagnon@ann.jussieu.fr\\ $^{2}$ GIREF, D\'epartement de math\'ematiques et de statistique, Universit\'e Laval, \\Qu\'ebec, Canada. jose.urquiza@mat.ulaval.ca }
\begin{document}

\maketitle

\begin{abstract}

For a Legendre-Galerkin semi-discretization of the 1-D homogeneous wave equation, the high frequency components  of the numerical solution prevent us from obtaining the boundary observability (inequality), uniformly with regard to the discretization parameter. A classical Fourier filtering that filters out the high frequencies is sufficient to recover the uniform observability. Unfortunately, this remedy needs to compute all the frequencies of the underlying system. In this paper we present three cheaper alternative remedies, namely a spectral filtering, a mixed formulation and Nitsche's method to append Dirichlet type boundary conditions. Our numerical results show indeed that uniform boundary observability inequalities may be recovered. On another hand, surprisingly, none of them seems to provide a uniform direct (or trace) inequality, a property which is needed in some existing general convergence results for the adjoint boundary controllability problem. We finally show in a numerical example that, despite this fact, convergence of the control approximations holds whenever the uniform observability inequality is observed numerically.
\end{abstract}

\section{Introduction}

Let us consider the one-dimensonal wave equation with homogeneous Dirichlet boundary conditions,  
\begin{equation}\label{obs}
\begin{cases}
u_{tt}(x,t)-u_{xx}(x,t)=0,  &  \quad \text{for } (x,t)\in(-1,1)\times (0,T)\\
u(-1,t)=u(1,t)=0,  & \quad \text{for } t\in(0,T)  \\
u(x,0)=u_0(x),u_t(x,0)=u_1(x), & \quad \text{for } x\in(-1,1)
\end{cases}
\end{equation}
With initial data $(u_0,u_1)\in H^1_0(-1,1)\times L^2(-1,1)$, 
the energy of the solution,
\begin{equation}\label{defenergie}
    E(u(t)):=\dfrac{1}{2}\int_{-1}^1 \left[ \left|u_t(x,t)\right|^2+\left|u_x(x,t)\right|^2 \right] \, \textrm{dx},
\end{equation}
is known to be conserved along time : $E(u(t))= E(u(0))$, $t\geq 0$. 

We are interested in the {\it boundary observability} of the wave equation. We say that it is {\it observable} in time $T>0$ if there exists $c_T>0$ such that 
\begin{equation}\label{intobs}
	c_T E(u(0)) \leq \int_0^T \left| u_x(1,t)\right|^2 \textrm{dt}, \quad \forall 		(u_0,u_1)\in H^1_0(-1,1)\times L^2(-1,1).
\end{equation}

Note that the right hand side of this inequality is finite for every time $T>0$. Indeed, for every $T>0$, there exists $C_T>0$ such that 
\begin{equation}\label{inttrace}
	\int_0^T \left| u_x(1,t)\right|^2 \textrm{dt} \leq C_TE(u(0)), \quad \forall 		(u_0,u_1)\in H^1_0(-1,1)\times L^2(-1,1).
\end{equation}
This is often referred as the {\it direct} or  {\it hidden regularity} inequality (\cite{JL88}). 

Whereas the direct inequality holds for every $T>0$, the observability inequality (\ref{intobs}) holds only if $T$ is sufficiently large, $T\geq 4$ in this case (see for instance \cite{KOMOR}). The time $T=4$ is the minimal time for all waves travelling within the spatial domain $(-1,\;1)$ at propagation speed 1 to reach the region of observation which is $x=1$. In fact, if $T=4$, inequalities \eqref{intobs} and \eqref{inttrace} hold with $c_T=C_T=8$, so that both inequalities become (the same) equalities.

Inequalities (\ref{intobs}) and (\ref{inttrace}) play a fundamental role in the Dirichlet boundary controllability of the wave equation 
\begin{equation}\label{cont} 
\begin{cases}
	y_{tt}(x,t)-y_{xx}(x,t)=0, & \quad \text{for } (x,t)\in(-1,1)\times (0,T), \\
	y(-1,t)=0 , &\quad   \text{for } t\in(0,T),  \\
	y(1,t)=v(t), &\quad \text{for } t\in(0,T) , \\
	y(x,0)=y_0(x),y_t(x,0)=y_1(x), & \quad \text{for } x\in(-1,1).
\end{cases}
\end{equation}
We say that (\ref{cont}) is {\it null controllable} in time $T>0$ if, for all initial data $(y_0,y_1)\in L^2(-1,1)\times H^{-1}(-1,1)$, there exists a control $v\in L^2(0,T)$ that drives the solution $y$ from $(y_0,y_1)$ to 
\begin{equation}\label{finaldata}
	y(.,T)=y_t(.,T)=0. 
\end{equation}
Because of the time-reversibility of the wave equation, null controllability is equivalent to controllability, which is the property to drive the solution of (\ref{cont}) from any initial state (at $t=0$) to any final state (at $t=T$) in $ L^2(-1,1)\times H^{-1}(-1,1)$. 

The equivalence between the observability of (\ref{obs}) and the controllability of (\ref{cont}) was shown in \cite{JL88} using the so-called Hilbert Uniqueness Method (HUM). Moreover, when controllable, a characterization of the control of minimal $L^2(0,T)$-norm was given. 

If $v$ is a control candidate, it must satisfy  
\begin{equation}\label{necess}
	\int_0^T v(t)\tilde{u}_x(1,t) \textrm{dt}=<(y_1,-y_0),(\tilde{u}_0,\tilde{u}_1)>, 	\quad \forall (\tilde{u}_0,\tilde{u}_1)\in H^1_0\times L^2,
\end{equation}
where $\tilde{u}$ is solution of \eqref{obs} with initial data $(\tilde{u}_0,\tilde{u}_1)$. To see this, one formally multiplies \eqref{cont} by $\tilde{u}(x,t)$ and make some integrations by parts. Note that here and in the sequel, $H^1_0$ stands for $H^1_0(-1,1)$, and the same notational simplification is adopted for spaces like $L^2(-1,1)$ and $H^{-1}(-1,1)$. Moreover, $<.,.>$ stands for the duality product between $H^{-1}\times L^2$ and $H^1_0\times L^2$.

The control of minimal $L^2(0,T)$-norm is simply given by 
\begin{eqnarray}\label{thecontrol}
	v(t)=u_x(1,t),
\end{eqnarray} 
where $u$ is the solution of (\ref{obs}) for suitably determined initial data  $(u_0,u_1)\in H^1_0\times L^2$.

How are $(u_0,u_1)$ determined? Define the linear, symmetric and positive operator
\begin{eqnarray}\label{defgram}
	\Lambda_T: H^1_0\times L^2 &\rightarrow & H^{-1}\times L^2\\
	<\Lambda_T(u_0,u_1),(\tilde{u}_0,\tilde{u}_1)> & = &\int_0^T u_x(1,t)		\tilde{u}_x(1,t) \textrm{dt}, \; \forall (\tilde{u}_0,\tilde{u}_1)\in H^1_0\times L^2.
\end{eqnarray} 
The direct inequality \eqref{inttrace} ensures that $\Lambda_T$ is well defined and continuous, while the observability inequality \eqref{intobs} ensures that $\Lambda_T$ is positive definite and thus invertible. Thus, if $T$ is sufficiently large for the observability inequality to hold, one can solve : find $(u_0,u_1)\in H^1_0\times L^2$ such that 
\begin{eqnarray}\label{eqgram}
	\Lambda_T(u_0,u_1)=(y_1,-y_0).
\end{eqnarray}

$\Lambda_T$ is referred as a {\it observability gramian} (for system \eqref{obs}) in the control literature (see for instance \cite{sontag}) or as the {\it HUM operator} (for control problem \eqref{cont}) in \cite{JL88}. Note also that, due to the symmetry of $\Lambda_T$, the solution  $(u_0,u_1)$ of \eqref{eqgram} can also be viewed as the minimizer over $H^1_0\times L^2$ of the quadratic functional 
\begin{equation}\label{qfunct}
	J(u_0,u_1)=\dfrac{1}{2}\int_0^T \left| u_x(1,t)\right|^2 \textrm{dt} \quad
	- \; <(y_1,-y_0),(u_0,u_1)>.
\end{equation}
From \eqref{eqgram} (or, equivalently, \eqref{qfunct}), one also deduces an upper bound for the control of the form 
$$
||v||_{L^2(0,T)}\leq \frac{C}{c_T}||y_0,y_1||_{L^2\times H^{-1}}.
$$

In the present work, we  consider space discretizations of the wave equation \eqref{obs} using a Legendre Galerkin method, with approximation vector spaces made of polynomials of degree $N$ or less. Let us denote by $u^N(x,t)$ the resulting approximation of $u$. 

We are interested with the discrete version of the observability inequality \eqref{intobs} that may be written :
\begin{eqnarray}\label{introintobsdisc}
	c_{N,T}E(u^N(0)) \leq \int_0^T \left| u^N_x(1,t)\right|^2 \textrm{dt} \; .
\end{eqnarray}
In particular, we would like the discrete observability inequality to hold uniformly with respect to the discretization parameter $N$, that is with 
$$\label{obsunif}
	c_{N,T}=c_{T}>0.
$$  

This uniform observability property is an important condition for proving the convergence (as $N\rightarrow\infty$) of the control approximations resulting from the approximation method chosen for \eqref{obs} (see for instance \cite{ZUA1,ErvedozaZuazua2012,ErvedozaZuazua2013}). An approximation $v_N$ of the control is obtained by minimizing the discrete analogue of the functional $J$, 
\begin{eqnarray}\label{mindisc}
	J_{N}(u_0^N,u_1^N)&=&\dfrac{1}{2}\int_0^T \left| u^N_x(1,t) \right|^2 \textrm{dt} \\ 
	&&-<(y_1,-y_0),(u_0^N,u_1^N)>_{H^{-1}\times L^2,H^1_0\times L^2}\nonumber,
\end{eqnarray}
over the approximation spaces for the initial data $(u_0^N,u_1^N)$, and then by setting $v_N(t)=u^N_x(1,t)$, 
where $u^N$ is the solution of \eqref{obsdiscret} with the initial data $(u_0^N,u_1^N)$ minimizing \eqref{mindisc}. A uniform observability property thus ensures that $J_N$ is uniformly  coercive and provides a uniform upper bound (with respect to the discretization parameter $N$) for $||v_N||_{L^2(0,T)}$.

Unfortunately, the usual convergence of solutions for approximations of wave equations like \eqref{obs} does not guarantee the convergence of the control approximations.  In particular, this does not guarantee the uniformity of \eqref{introintobsdisc}.  This was first proved in the context of finite-difference and finite element approximations in \cite{infantezuazua}.
Similarly, the standard Legendre Galerkin semi-discretization does not provide a uniform observability inequality (\cite{JU}). This pathology is caused by the high frequency components of the solutions,  the highest numerical eigenmodes not being suitable approximations of the corresponding eigenmodes of the continuous problem.  We refer the reader to \cite{ZUA1, ErvedozaZuazua2012} for interpretations in terms of numerical wave propagation and for a sharp analysis of the phenomena involved in the context of finite-difference and finite-element approximations of the one-dimensional (and two-dimensional) wave equation.

For all these standard space discretization methods, a remedy to recover the uniformity of the observability is to filter out the high frequency components of the solution. This method is often referred as {\it Fourier filtering}. It was shown in \cite{infantezuazua,ZuazuaJMAP,JU} that this method led to the uniformity of the observability inequality, for the numerical scheme considered in these papers. This method, however, has a major drawback : it is necessary to compute the eigenmodes (eigenvalues and associated eigenfunctions), which can be costly, especially if one wants to consider 2 or 3-dimensional spatial domains.  

Based on numerical experiments in the case of a finite-element or a finite-difference semi-discretization of (\ref{obs}), other remedies were suggested by Glowinski in \cite{GLL,GLO}. Among them, let us cite a Tychonoff regularization, a mixed formulation of (\ref{obs}), a bi-grid method and a penalty method. It was shown in \cite{NEGZ} that a bi-grid method applied to a finite-difference semi-discretization allows one to obtain the uniformity of the observability inequality and the convergence of the numerical controls. In \cite{MICU}, the same positive results were proved using a mixed finite element approximation. 

In this paper we propose three remedies to recover the uniformity of the boundary observability inequality in the case of a spectral semi-discretization of (\ref{obs}). In the context of the Legendre-Galerkin semi-discretization of (\ref{obs}) considered in \cite{JU}, we show numerically that a spectral filtering (a commonly used procedure in spectral methods), a mixed method and  Nitsche's method to append Dirichlet boundary conditions allow one to recover the uniformity of $c_{N,T}$. To our knowledge, this is the first time that one obtains the uniformity of $c_{N,T}$ without the use of Fourier filtering in the case of a spectral semi-discretization of (\ref{obs}). 


Unfortunately, these remedies are not fully satisfactory since they do not seem to guarantee  the  uniformity of the discrete analogue of the direct inequality \eqref{inttrace}, 
\begin{eqnarray}
\int_0^T \left| u^N_x(1,t)\right|^2 \textrm{dt} \leq C_{N,T}E(u^N(0)), \label{introinttracedisc} 
\end{eqnarray}
that is with 
$$
	C_{N,T}<C_T,
$$
independently of the discretization parameter $N$.  In this regard, the situation here is very different from what happens with classical finite element or finite difference methods, where the uniform direct inequality holds even without filtering out the highest frequencies (\cite{infantezuazua}). 

These negative results, if confirmed theoretically, would prevent us from using the abstract convergence results (for the control and the controlled solution)
 in \cite{ErvedozaZuazua2012,ErvedozaZuazua2013} which need both the uniform observability and the uniform direct inequality as well. While the uniform observability property seems essential (see \cite{ErvedozaZuazua2012}, section 4.4.) the necessity of the uniform direct inequality for the convergence of the control approximations is not clear to us. 
 
Nevertheless, we will show a numerical example where convergence does not hold with the classical Legendre Galerkin method (without uniform observability), whereas convergence does hold if we compute the approximate control using one of the remedies (new discrete formulations with uniform observability) proposed in this paper. 

The paper is organized as follows. In Section 2, we recall some results on the behaviour of $c_{N,T}$ and $C_{N,T}$ obtained in \cite{JU} with the Legendre-Galerkin approximation of (\ref{obs}). In Section 3, we show numerically that a modification of the discrete version of the functional $J(u_0,u_1)$ via a spectral filtering allows us to obtain a uniform and positive lower-bound on $c_{N,T}$. In Section 4, we obtain numerically the uniformity of $c_{N,T}$ when considering a mixed Legendre-Galerkin method. The uniform result on $c_{N,T}$ differs from the one obtained in \cite{MICU} since a weaker version of the observability inequality was necessary. In Section 5, we introduce Nitsche's Method to append Dirichlet type boundary conditions. Two variants of Nitsche's method are used for both the controlled wave equation and the dual observability problem, and again the associated observability property is shown to hold uniformly in our numerical experiments. Finally, in Section 6, we present a numerical example of problem to be controlled, showing the convergence of the control approximations for all the remedies proposed in this paper.  

\section[Legendre Galerkin semi-discretization]{Legendre Galerkin approximation and discrete observability pro\-perties}

\subsection{Legendre Galerkin approximation}

First of all, let us recall that the solution of the wave equation admits the Fourier expansion 
\begin{equation*}
 	u(x,t)=\sum_{k=1}^{\infty} \left(\alpha_{k}\cos(\sqrt{\lambda_{k}}t)+ 					\dfrac{\beta_{N}}{\sqrt{\lambda_{k}}} \sin(\sqrt{\lambda_{k}}t)\right) \phi_{k}(x),
\end{equation*}
where 
\begin{eqnarray*}
	\alpha_{k}=(u_0,\phi_{k})_{L^2(-1,1)}; \quad 
	\beta_{k}=(u_1,\phi_{k})_{L^2(-1,1)},
\end{eqnarray*}
and $(\lambda_k,\phi_k)_{k\geq 1}$ is the eigensystem of the eigenproblem : 
 find $\lambda\in\mathbb{R}$, $\phi\in H^1_0$ such that 
\begin{equation}\label{probproprecontinu}
	(\phi,\psi)_{H^1_0(-1,1)}=\lambda(\phi,\psi)_{L^2(-1,1)}, 
	\quad \forall \psi\in H^1_0(-1,1).
\end{equation}
with $(\phi_k,\phi_k)_{L^2(-1,1)}=1$, that are 
\begin{eqnarray*}
	\lambda_{k}=\left(k\pi\right)^2 \quad 
	\phi_{k}(x)=\sin(\sqrt{\lambda_k}\,x).
\end{eqnarray*}

Let us now briefly describe the semi-discretization of (\ref{obs}) studied in \cite{JU}.
Let $L_k$ denote the $k$-th order Legendre polynomial, given by the recurrence relation,
\begin{equation*}
\begin{array}{l} 
	L_0(x)=1, \\
	L_1(x)=x, \\
	(k+1)L_{k+1}(x)=(2k+1)xL_{k}(x)-kL_{k-1}(x), \quad k\geq 2.
\end{array}
\end{equation*}
Let $\PP_N^0(-1,1)$ denote the space of polynomials of degree at most $N$ vanishing at $x=-1$ and $x=1$. We notably have (see \cite{SHEN}) : 
\begin{eqnarray*}
	\PP_N^0(-1,1)=\textrm{span}\{\widetilde{L}_1(x),\ldots,\widetilde{L}_{N-1}(x)\},
\end{eqnarray*} 	
where
\begin{equation*}
	\widetilde{L}_{k}(x)=c_k(1-x^2)L_k'(x), \quad c_k=\dfrac{\sqrt{k+1/2}}{k(k+1)}, \quad k\geq 1.
\end{equation*}

The Legendre-Galerkin semi-discretization of (\ref{obs}) is then  
\begin{eqnarray}\label{obsdiscret}
	\int_{-1}^1 u^N_{tt}(x,t)\psi^N(x) \, \textrm{dx}+\int_{-1}^1 u^N_x(x,t)\psi^N_x(x)\, \textrm{dx} = 0, \; \forall \psi^N		\in \PP_N^0(-1,1), \; t\in(0,T)
\end{eqnarray}
with suitable initial conditions. Here, the approximation $u^N(.,t)$ of $u(.,t)$  is sought in $\PP_N^0(-1,1)$, that is in the form 
\begin{equation*}
 u^N(x,t)=\sum_{k=1}^{N-1} a_{N,k}(t)\widetilde{L}_k(x)
\end{equation*}
Introducing $\textbf{a}_N(t):=(a_{N,1}(t),\ldots,a_{N,N-1}(t))^t$, system (\ref{obsdiscret}) can be written in the matrix form
\begin{equation*}
M_N\textbf{a}_N''(t)+K_N\textbf{a}_N(t)=0, 
\end{equation*}
where $M_N$ and $K_N$ are positive definite symmetric matrices given by (see \cite{SHEN})
\begin{equation*}
\begin{array}{l}
M_N(i,j)=\left\{\begin{array}{l l}\frac{2}{(2i-1)(2i+3)} & \textrm{ if i=j}\\
         \frac{-1}{(2r+3)\sqrt{(2r+1)(2r+5)}} & \textrm{ if } |i-j|=2\\
	 0 & \textrm{ otherwise } 
         \end{array}\right.\\
K_N(i,j)=\left\{\begin{array}{l l} 1 & \textrm{ if } i=j\\
	 0 & \textrm{ otherwise,}
         \end{array}\right.
\end{array}
\end{equation*}
with $r=\min\{i,j\}$. 

Taking $\psi^N=u_t^N(x,t)$ in (\ref{obsdiscret}), one obtains that the energy of the solution of (\ref{obsdiscret})
is preserved along time : $E(u^N(t))=E(u^N(0)), \forall t\geq 0$. 

\begin{figure}[!th]
\begin{center}
	\includegraphics[height=8cm]{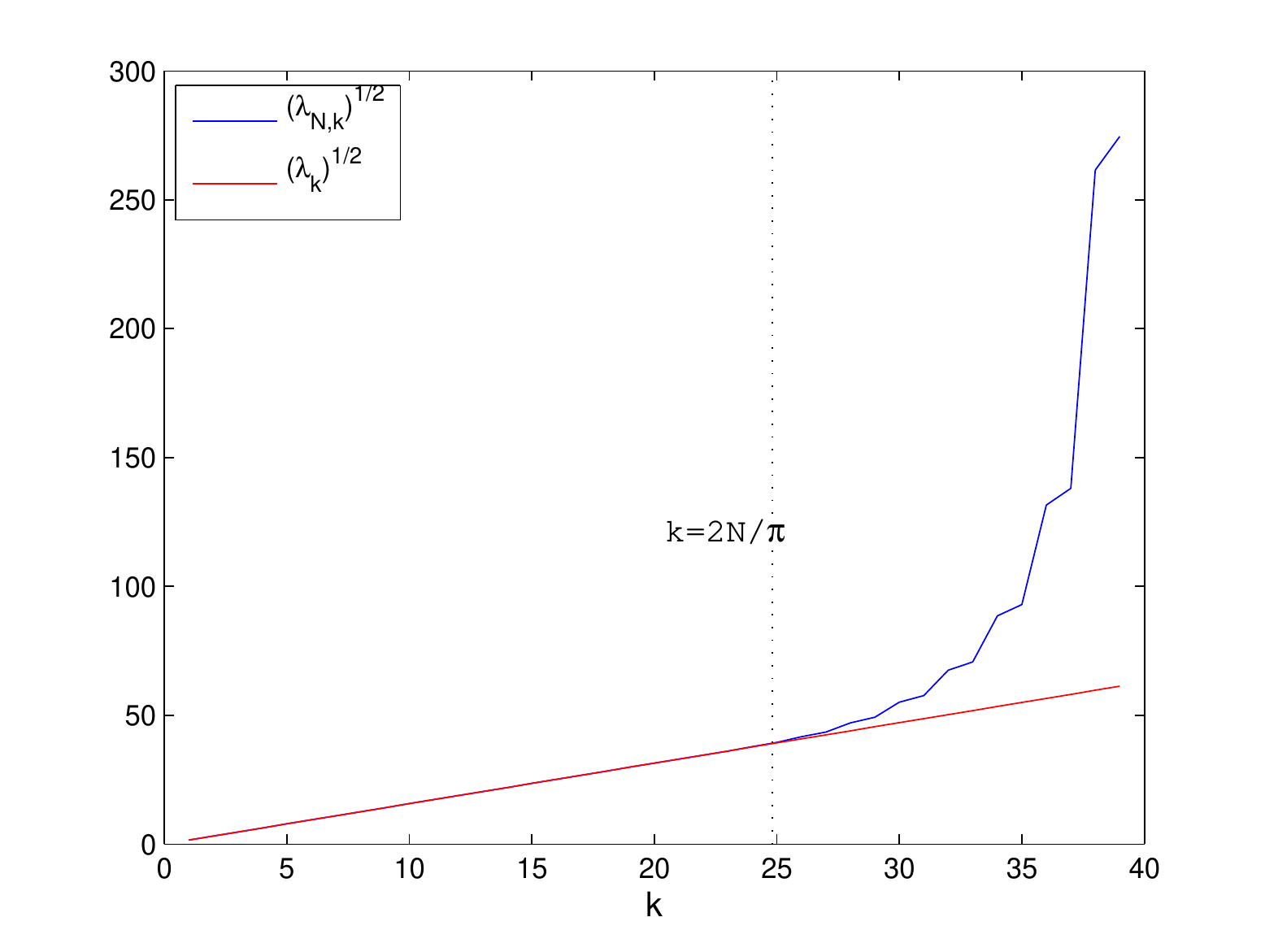}
\end{center}
\caption{\footnotesize{Behaviour of the square roots of the discrete eigenvalues and of the continuous ones for $N=40$.}}\label{vapeqdiscret}
\end{figure}

Let $(\lambda_{N,k},\phi_{N,k})$, $1\leq k\leq N-1$ denote the eigenpair solutions of the eigenproblem: find $\lambda\in\mathbb{R}$, $\phi\in \PP_N^0(-1,1)$ such that 
\begin{equation}\label{probproprediscret}
(\phi,\psi^N)_{H^1_0(-1,1)}=\lambda(\phi,\psi^N)_{L^2(-1,1)}, \quad \forall \psi\in \PP_N^0(-1,1).
\end{equation}
It was shown in \cite{JU} that the eigenvalues $\{\lambda_{N,k}\}$ of the discrete eigenproblem are the zeros of  
\begin{equation*}
\dfrac{J_{\nu}(x)}{J_{-\nu}(x)}=(-1)^N\tan\left(x-\dfrac{N\pi}{2}\right), 
\end{equation*}
where $\nu=N+1/2$ and where $J_{\nu}$ stands for the spherical Bessel's function of order $\nu$. 
Figure \ref{vapeqdiscret} illustrates the behaviour of $\sqrt{\lambda_{N,k}}$, for $N=40$, in comparison with $\sqrt{\lambda_{k}}=k\pi$. The first $2/\pi$ fraction of them are very well approximated, as proved in \cite{VANDEVEN}.

The Fourier expansion of the solution of (\ref{obsdiscret}) on the orthonormal basis formed by the discrete normalized (for the $L^2(-1,1)$ norm) eigenfunctions $\phi_{N,k}$ writes 
\begin{equation*}
 	u^N(x,t)=\sum_{k=1}^{N-1} \left(\alpha_{N,k}\cos(\sqrt{\lambda_{N,k}}t)+ 				\dfrac{\beta_{N,k}}{\sqrt{\lambda_{N,k}}} \sin(\sqrt{\lambda_{N,k}}t)\right)				\phi_{N,k}(x),
\end{equation*}
where, for $1\leq k \leq N-1$,
\begin{eqnarray*}
	\alpha_{N,k}=(u_0^N,\phi_{N,k})_{L^2(-1,1)}, \quad \beta_{N,k}=(u_1^N,			\phi_{N,k})_{L^2(-1,1)}.
\end{eqnarray*}
We also consider its Fourier truncation of order $M$, $1\leq M \leq N-1$,
\begin{equation*}
 	u^{N,M}(x,t)=\sum_{k=1}^{M} \left(\alpha_{N,k}\cos(\sqrt{\lambda_{N,k}}t)+ 		\dfrac{\beta_{N,k}}{\sqrt{\lambda_{N,k}}} \sin(\sqrt{\lambda_{N,k}}t)\right)		\phi_{N,k}(x).
\end{equation*}

\subsection{Discrete observability}

\begin{figure}[!th]
\begin{center}
	\includegraphics[height=4cm]{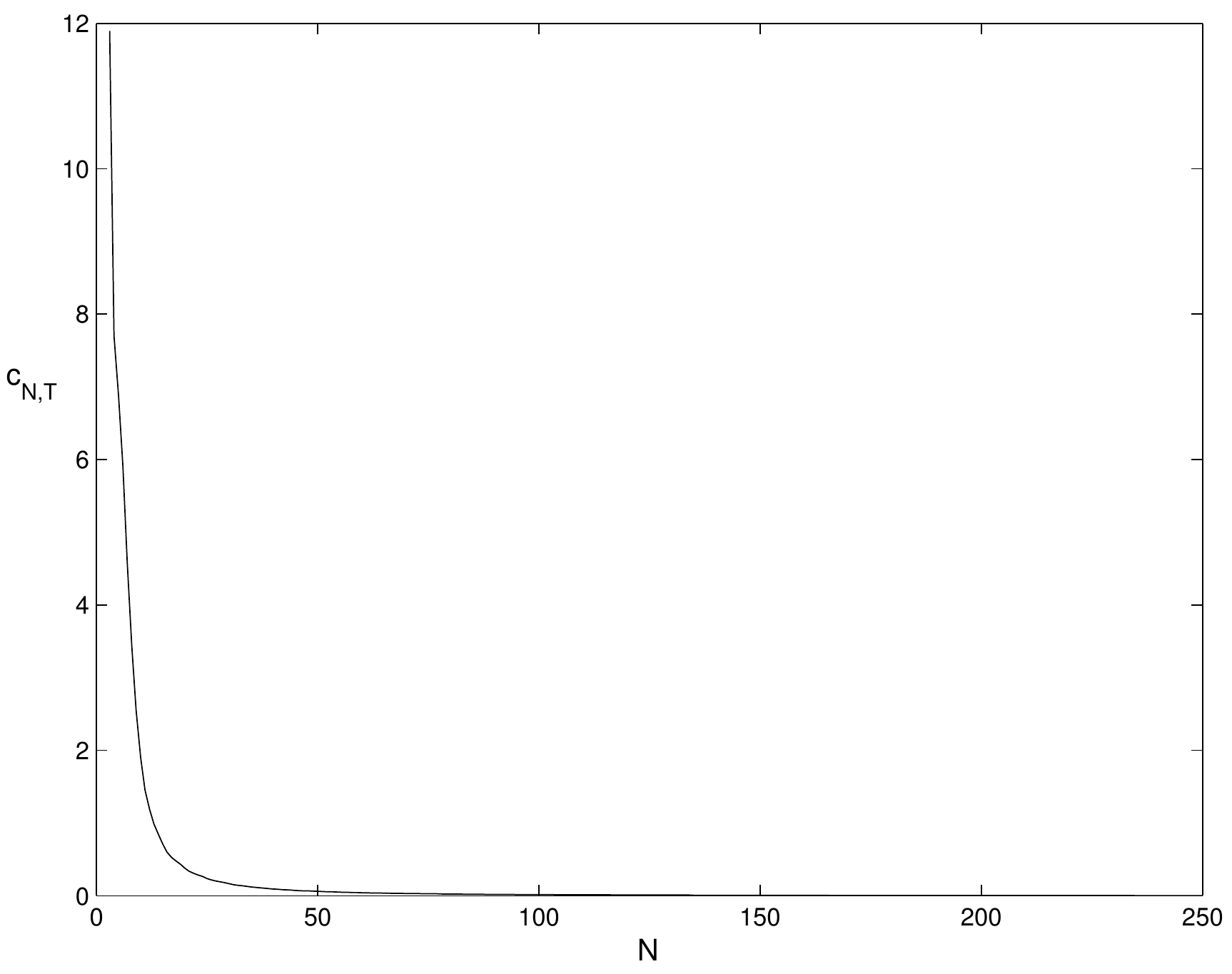} \quad \includegraphics[height=4cm]{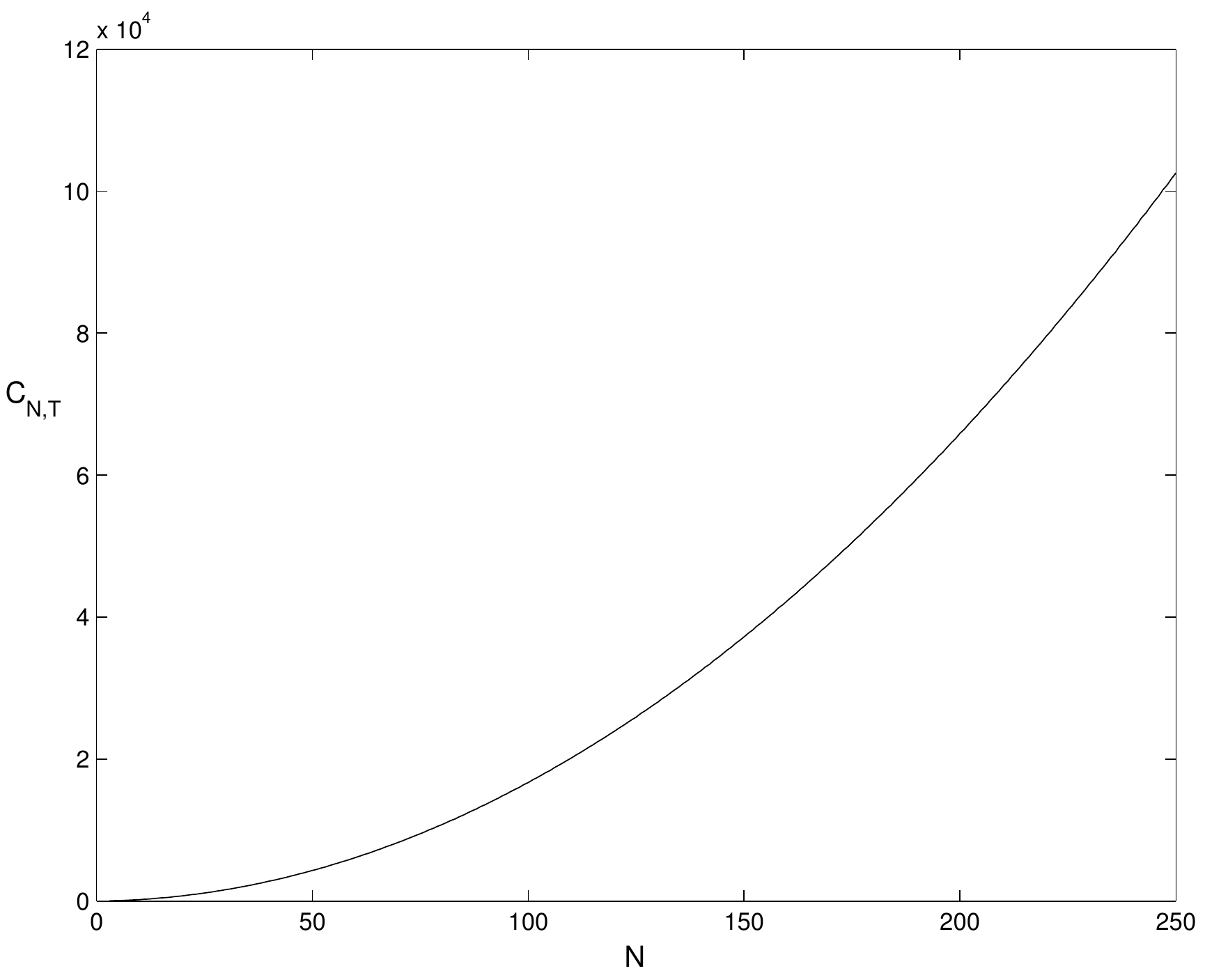} 
\end{center}
\caption{\footnotesize{Behaviour of $c_{N,T}$ (left) and $C_{N,T}$ (right) with $T=8$}}\label{conststandard}
\end{figure}

Considering the discrete observability and direct inequalities \eqref{intobs} and \eqref{inttrace}, we redefine the constant $c_{N,T}$ and $C_{N,T}$ as the best possible ones :
\begin{equation}\label{cnt1}
c_{N,T} = \inf_{(u_0^N,u_1^N)\in {(\PP_N^0)}^2} \left(\int_0^T \left| u^N_x(1,t)\right|^2 dt \right) \;/ \; E(u^N(0)),
\end{equation}
\begin{equation}\label{cnt2}
C_{N,T} = \sup_{(u_0^N,u_1^N)\in {(\PP_N^0)}^2} \left(\int_0^T \left| u^N_x(1,t)\right|^2 dt \right) \;/ \; E(u^N(0)),
\end{equation}
As mentioned in the introduction, the observability and direct inequalities do not hold uniformly, as we have for every $T>0$ (\cite{JU}): 
%

\begin{equation*}
 c_{N,T} \rightarrow 0 \quad \textrm{ and } \quad C_{N,T}\rightarrow \infty \quad \textrm{ as } \quad N\rightarrow \infty. 
\end{equation*}
Figure \ref{conststandard} shows the behaviour of $c_{N,T}$ and $C_{N,T}$ with respect to  $N$ for $T=8$. We refer to the Appendix for the numerical method used to compute $c_{N,T}$ and $C_{N,T}$.

Nevertheless, if we consider a suitable Fourier truncation $u^N_M$, which filters out the high frequency components, then we can recover uniform observability and direct inequalities. More precisely, if $T>4$ and if $M<\frac{2}{\pi}N$, that is if we keep only the $2/\pi$ fraction of the lowest frequency components, then 
\begin{equation*}
	c_T \;E(u^{N,M}(0))\,\leq \,\int_0^T \left| \dfrac{\partial u^{N,M}}{\partial x}(1,t)\right|	^2 dt \,\leq \, C_T \;E(u{N,M}(0))
\end{equation*}
where $c_T>0$ and $C_T>0$ are independent of $N$ (\cite{JU}). 

\begin{figure}[!th]
\begin{center}
	\includegraphics[height=6cm]{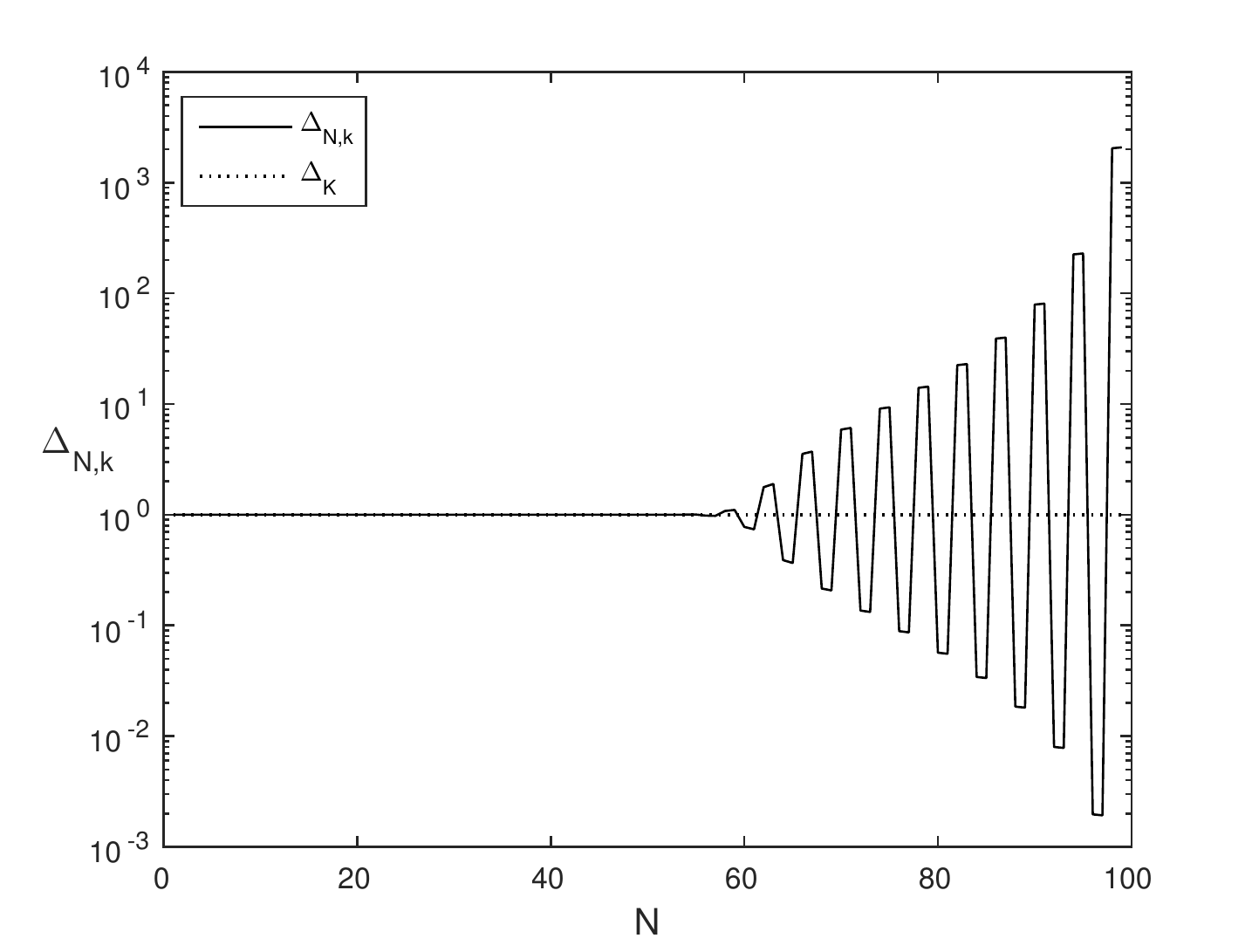}
\end{center}
\caption{\footnotesize{Behaviour of $\Delta_{N,k}$, for $N=40$, compared to their continuous counterparts  $\Delta_{k}=1$}}\label{figdelta}
\end{figure}

The proof of this result relies on two ingredients: First, the existence of a uniform gap between consecutive eigenvalue square roots, 
\begin{equation}\label{gapuniforme}
\sqrt{\lambda_{N,k}} - \sqrt{\lambda_{N,k-1}}\geq \delta>0
\end{equation}
which is proved in \cite{JU} to hold as long as $k<\frac{2}{\pi} N$; secondly, on the uniform lower and upper bounds on
\begin{equation*}
\Delta_{N,k}=\dfrac{\left| \phi_{N,k}'(1)\right|^2 }{\displaystyle \int_{-1}^1\left|\phi_{N,k}'(x)\right|^2\textrm{ dx}},
\end{equation*}
that is 
\begin{equation}\label{Deltauniforme}
c\leq\Delta_{N,k} \leq C, 
\end{equation}
which also holds as long as $k<\frac{2}{\pi} N$. 
Note that $\Delta_{N,k}$ has to be compared with its counterpart in the continuous case, which is 
\begin{equation*}
\Delta_{k}=\dfrac{\left| \phi_{k}'(1)\right|^2 }{\displaystyle \int_{-1}^1\left|\phi_{k}'(x)\right|^2\textrm{ dx}}.
\end{equation*}
A straightforward computation shows that $\Delta_{k}=1$, for every $k\geq 1$.  
Figure \ref{figdelta} illustrates the fact that  values of $\Delta_{N,k}$ for $N=40$ are indeed close to the exact values $\Delta_{k}=1$ as long as $k<\frac{2}{\pi} N$; 

If no filtering is applied, while the uniform gap property \eqref{gapuniforme}
between consecutive eigenvalue square roots seems to hold still (see Figure \ref{vapeqdiscret}), Figure \ref{figdelta} shows that both the lower and upper uniform bound property  on $\Delta_{N,k}$ in \eqref{Deltauniforme} do not hold anymore, as proved in \cite{JU}.

Let us point that the situation differs when using finite difference or finite element methods. In these cases, there is no uniform lower gap between consecutive eigenvalue square roots (caused by the highest numerical eigenvalues) causing by itself the non-uniform observability inequality, while a uniform upper bound on the equivalent $\Delta_{N,k}$ alone ensures a uniform direct (continuity) inequality.

\section{Spectral Filtering}

We have just recalled that the lack of uniform observability is due to the high frequency components of the approximate solution of the wave equation. On another hand, spectral approximation methods, like the Legendre Galerkin method, are known to be very precise when exact solutions are very smooth but also to lose part of their efficiency when solutions are not that smooth. In the latter case, the high frequency components of the approximation are again seen as the main source of errors or instabilities, for instance giving rise to the famous Gibbs phenomenon for discontinuous exact solutions.   

A way to recover precision is to dissipate the high modes using {\it spectral filters} (see for instance \cite{CAN1,HEST}). These filters are called {\it spectral filters} because the filtered approximate solutions are obtained by modifying their spectral (here, polynomial) expansion. Given a spectral expansion 
 \[
u^N(x)=\sum_{k=1}^{N-1} \alpha_{N,k}\widetilde L_k(x),  
\] 
the filtered expansion $\FF_N u^N$ is 
\begin{equation*}
\FF_N u^N(x)=\sum_{k=1}^{N-1} \sigma\left(\dfrac{k-1}{N-1}\right) \alpha_{N,k}\widetilde{L}_k(x)
\end{equation*}
where $\sigma$ is called a {\it filter function} or simply a {\it filter}. $\sigma(\eta)$ has values  close to $1$ for small values of $\eta$ and close to $0$ for $\eta$ approaching $1$.

More precisely, we say that 
$\sigma : [0,\,1]\rightarrow [0,\,1]$ is a $p${\it-th order filter} ($p\geq 1$) if
 \begin{enumerate}
	 \item [(i)] $\sigma \in C^{p-1}([0,\, 1])$,
 	\item [(ii)] $\sigma(0)=1$, $\sigma^{(j)}(0)=0$ for $1\leq j \leq p-1$,
  	\item [(iii)] $\sigma(1)=0$, $\sigma^{(j)}(1)=0$ for $1\leq j \leq p-1$. 
\end{enumerate}

Legendre polynomials are the eigenfunctions of a Sturm-Liouville operator arranged in the increasing order of their associated eigenvalues or frequencies. The idea behind this filtering procedure is thus to keep almost unchanged the low frequency components of the expansion while having a progressive damping on the higher modes. 

Note that the filtering is generally performed on the expansion in terms of the Legendre polynomials $L_k(x)$, $k=0,...,N$, which are of increasing frequency with increasing $k$, instead of $\widetilde{L}_{k}(x)$. But as we have 
\begin{eqnarray}
\widetilde{L}_{k}(x)=\frac{1}{\sqrt{4k+2}} (L_{k-1}(x)-L_{k+1}(x)), \quad k\geq 1,
\end{eqnarray}
functions $\widetilde{L}_{k}(x)$ are still of increasing frequency components with increasing $k$. 

In this work we consider several filter functions, all presented and discussed in \cite{CAN1}:
\begin{itemize}
\item \textbf{Ces\'aro filter} (first order)
\begin{equation*}
	\sigma(\eta)=1-\eta, 
\end{equation*}

\item \textbf{Lanczos filter} (first order)
\begin{equation*}
	\sigma(\eta)=\dfrac{\sin\left(\pi \eta\right)}{\pi \eta}, 
\end{equation*}

\item \textbf{Raised cosine filter} (second order)
\begin{equation*}
	\sigma(\eta)=\dfrac{1+\cos\left(\pi \eta\right)}{2}, 
\end{equation*}

\item \textbf{Sharpened raised cosine filter} (8-th order)
\begin{equation*}
	\sigma(\eta)=\sigma_0\left(\eta\right)^4\left(35-84\sigma_0\left(\eta\right)			 	+70\sigma_0\left(\eta\right)^2-20\sigma_0\left(\eta\right)^3 \right), 
\end{equation*}
where $\sigma_0$ denotes the raised cosine filter. 

\item \textbf{Vandeven filter} (p-th order)
\begin{equation*}
	\sigma(\eta)=1-\dfrac{(2p-1)!}{(p-1)!^2}\int_0^\eta (t(1-t))^{p-1} \, \textrm{dt} 
\end{equation*}

\item \textbf{Exponential filter} (p-th order)
\begin{equation*}
	\sigma(\eta)=e^{-\alpha \eta^p},  \quad \quad \alpha>0.
\end{equation*}
\end{itemize}

The exponential filter does not satisfy condition (iii) of the definition. However, taking $\alpha=-\log(\epsilon_m)$, where $\epsilon_m$ is the machine accuracy, the value of $\sigma(1)$ will be computationally interpreted as a 0. The exponential filter has the advantage to be of arbitrary order, like the Vandeven filter, but compared to the latter it has a much lower computational cost. 

Vandeven and exponentiel filters in the context of Legendre spectral methods were studied numerically in \cite{HEST}. There, it was observed that they behave similarly, they improve accuracy in polynomial expansions of smooth functions and they may be essential to stability for first order wave equations. It was also observed that increasing the order of the filter improves its accuracy. 

\begin{figure}[!th]
\title{}
\begin{center}
\begin{tabular}{c}
 \mbox{\includegraphics[height=6cm]{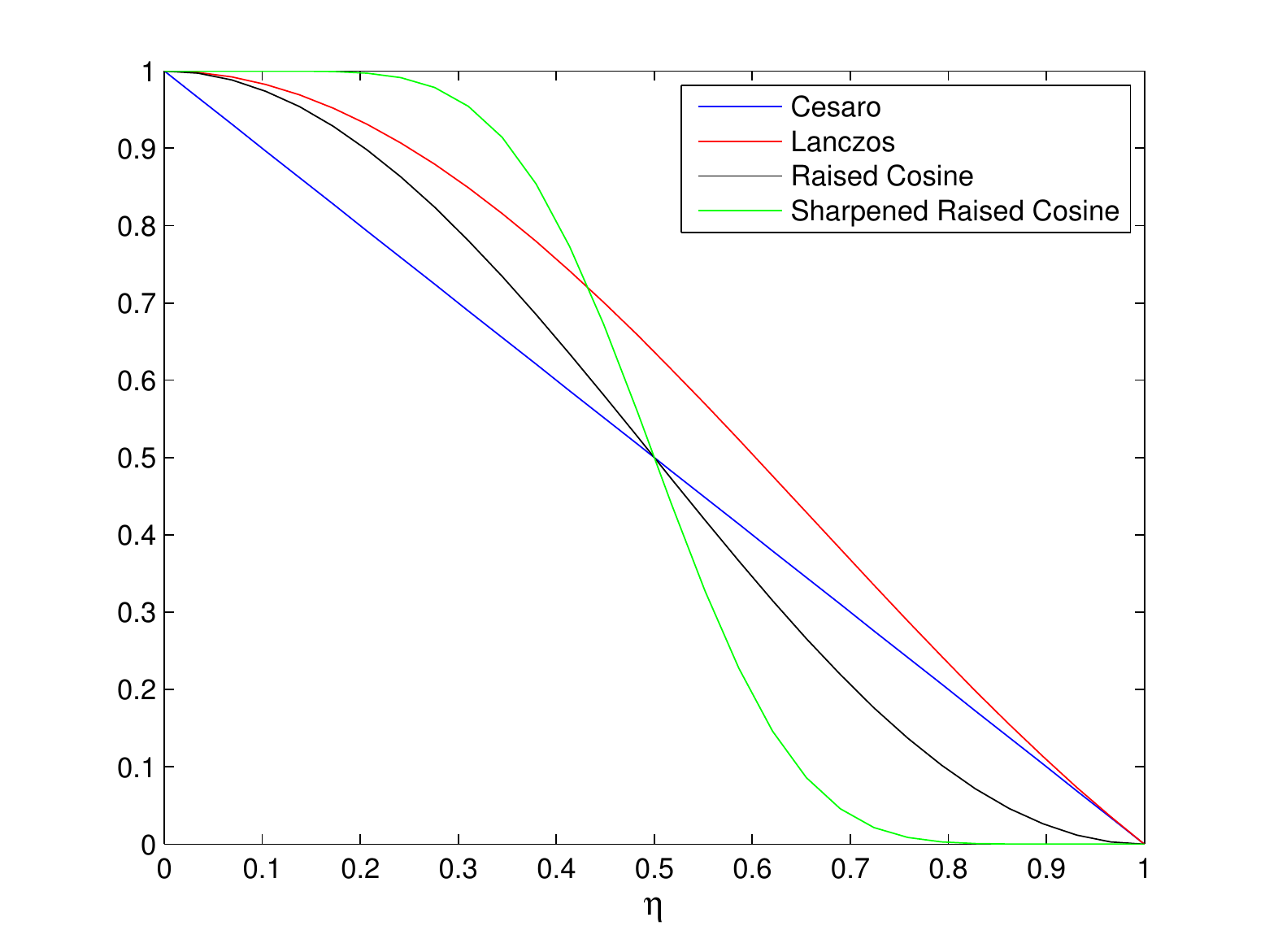}} 
\end{tabular}
\end{center}
\caption{\footnotesize{Graphs of Ces\`aro, Lanczos, raised cosine and sharpened raised cosine filters}}\label{filters}
\end{figure}

\begin{figure}[!th]
\title{}
\begin{center}
\begin{tabular}{c}
 \mbox{\includegraphics[height=5cm]{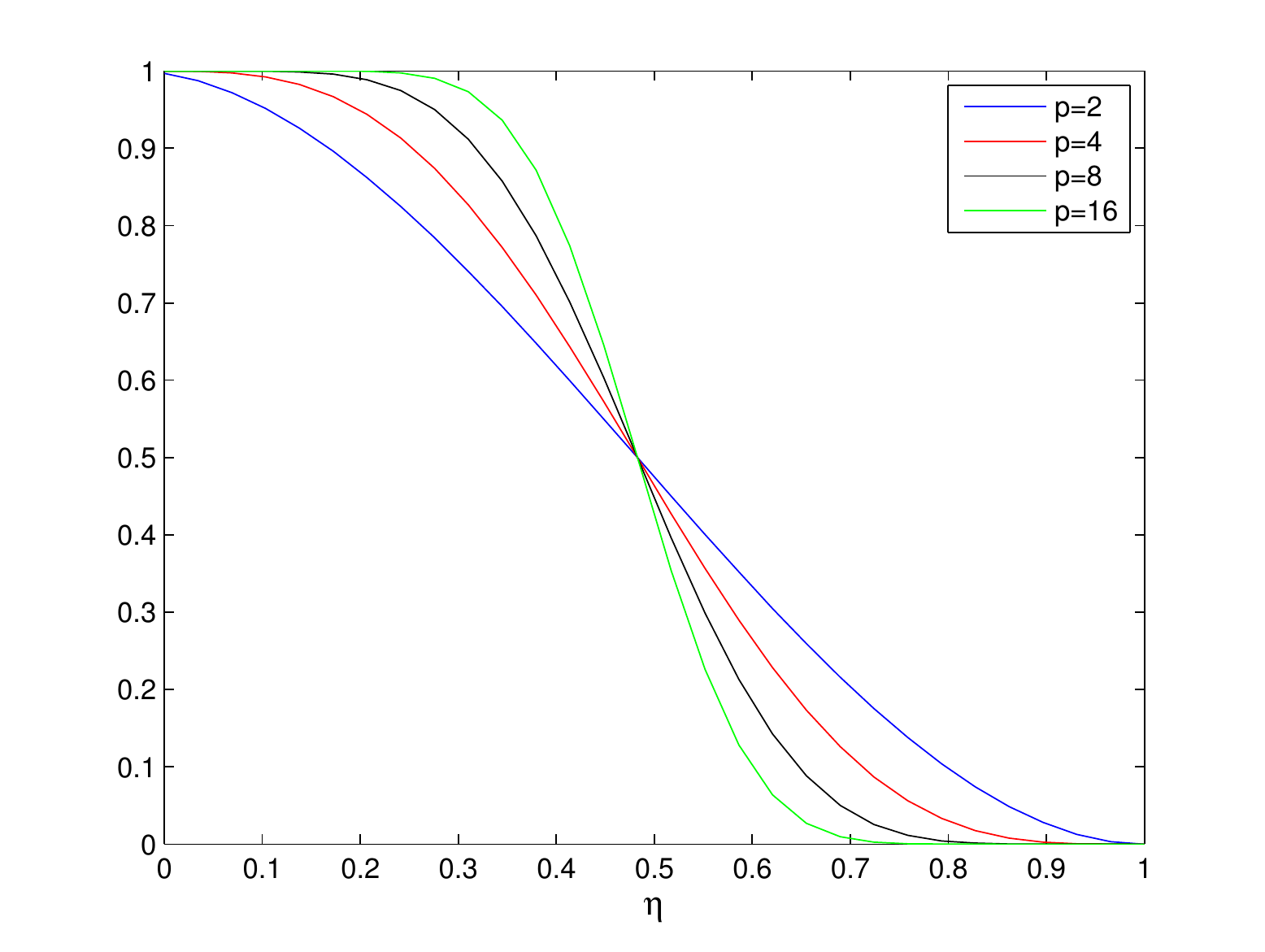}} \
 \mbox{\includegraphics[height=5cm]{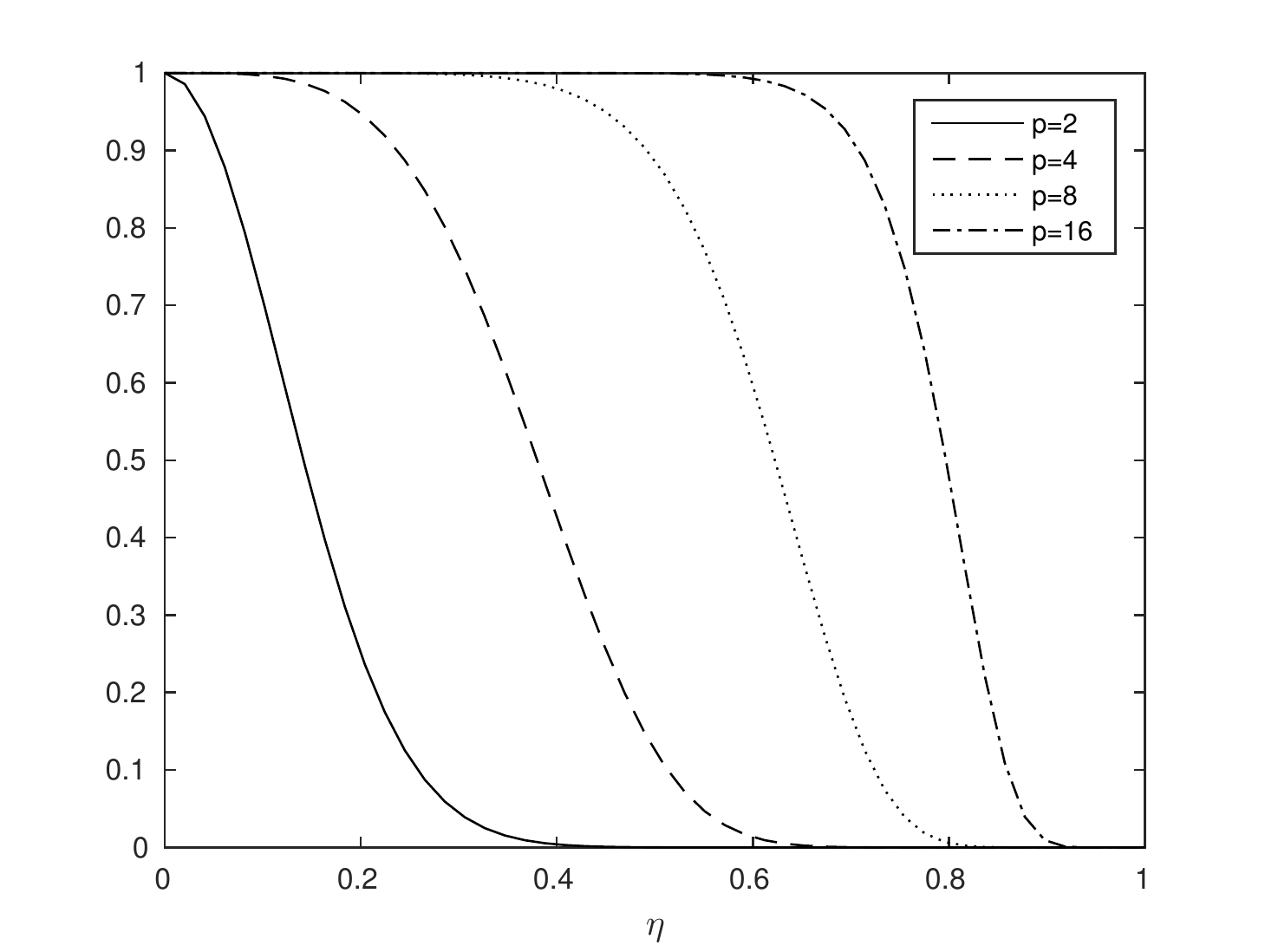}}
\end{tabular}
\end{center}
\caption{\footnotesize{Graphs of Vandeven (left) and exponential (right) filters for different values of $p$}}\label{filtersp}
\end{figure}

Figure \ref{filters} shows the graphs of the first four filters, whereas Figure \ref{filtersp} shows the graphs of the last two, for different values of $p$.

Regarding our observability problem, we now consider the boundary observation of the filtered expansion $\FF_N u^N$ of $u^N$, that is 
\begin{equation*}
(\FF_N u^N)_x(1,t)=\sum_{k=1}^{N-1} \sigma\left(\dfrac{k-1}{N-1}\right) a_{N,k}(t)\widetilde{L}_k'(1).
\end{equation*}
$(\FF_N u^N)_x(1,t)$ is expected to converge toward $u_x(1,t)$ (for an appropriate norm). 
The discrete controllability problem that we consider reduces to minimizing the functional 
\begin{eqnarray}\label{minfiltre}
	J_{N,\sigma}(u_0^N,u_1^N)&=&\dfrac{1}{2}\int_0^T \left| (\FF_Nu^N)_x(1,t) \right|^2 \textrm{dt} \\ 
	&&-<(y_1,-y_0),(u_0^N,u_1^N)>_{H^{-1}\times L^2,H^1_0\times L^2}\nonumber
\end{eqnarray}
over $(\PP_N^0)^2$. An approximation of the exact control then is $v_N(t)=(\FF_N u^N)_x(1,t)$, where $u^N$ is the solution of \eqref{obsdiscret} with the initial data $(u_0^N,u_1^N)\in(\PP_N^0)^2$ minimizing \eqref{minfiltre}.
 

We now look for uniform  bounds (with respect to $N$) of the best constants $c_{N,T}$ and $C_{N,T}$ appearing in
\begin{eqnarray}
 c_{N,T}E(u^N(0)) \leq \int_0^T \left| \left(\FF_Nu^N\right)_x(1,t)\right|^2 dt \label{obsfiltre} \\
\int_0^T \left| \left(\FF_Nu^N\right)_x(1,t)\right|^2 dt \leq  C_{N,T}E(u^N(0)). \label{tracefiltre}
\end{eqnarray}

\begin{figure}[!th]
\title{}
\begin{center}
\begin{tabular}{c}
 \mbox{\includegraphics[height=5cm]{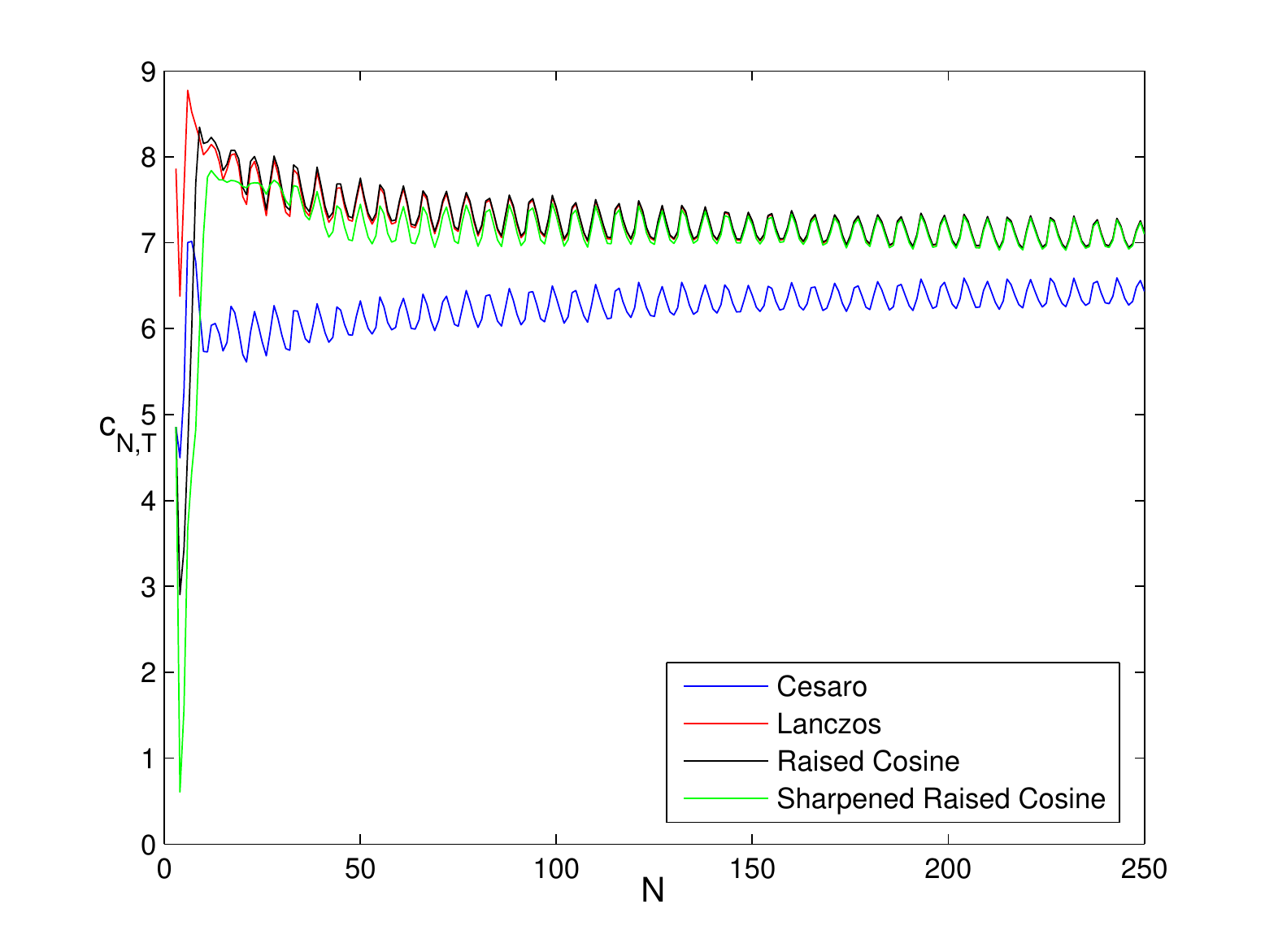}} 
 \mbox{\includegraphics[height=5cm]{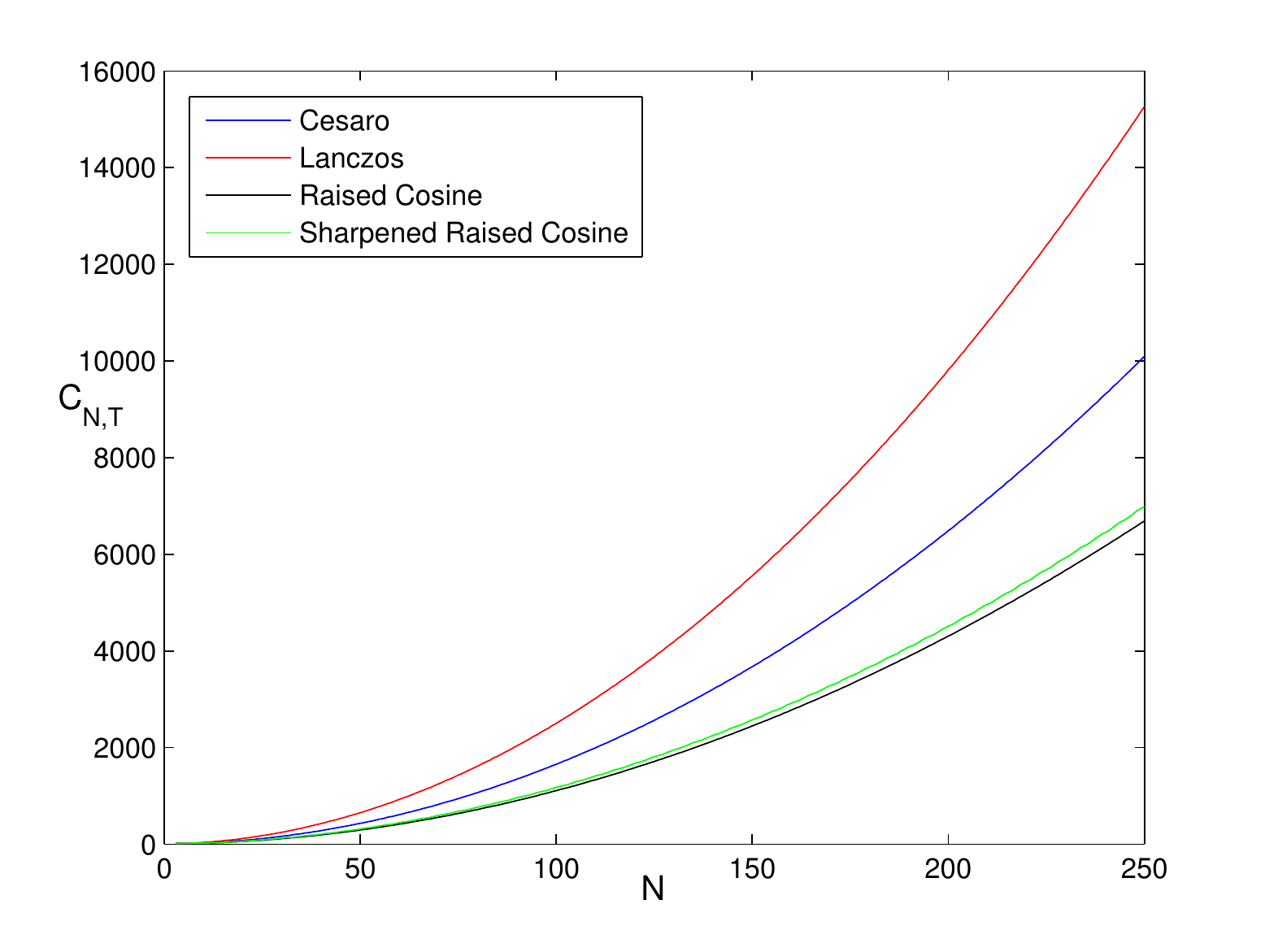}} \\
\end{tabular}
\begin{tabular}{c}
 \mbox{\includegraphics[height=5cm]{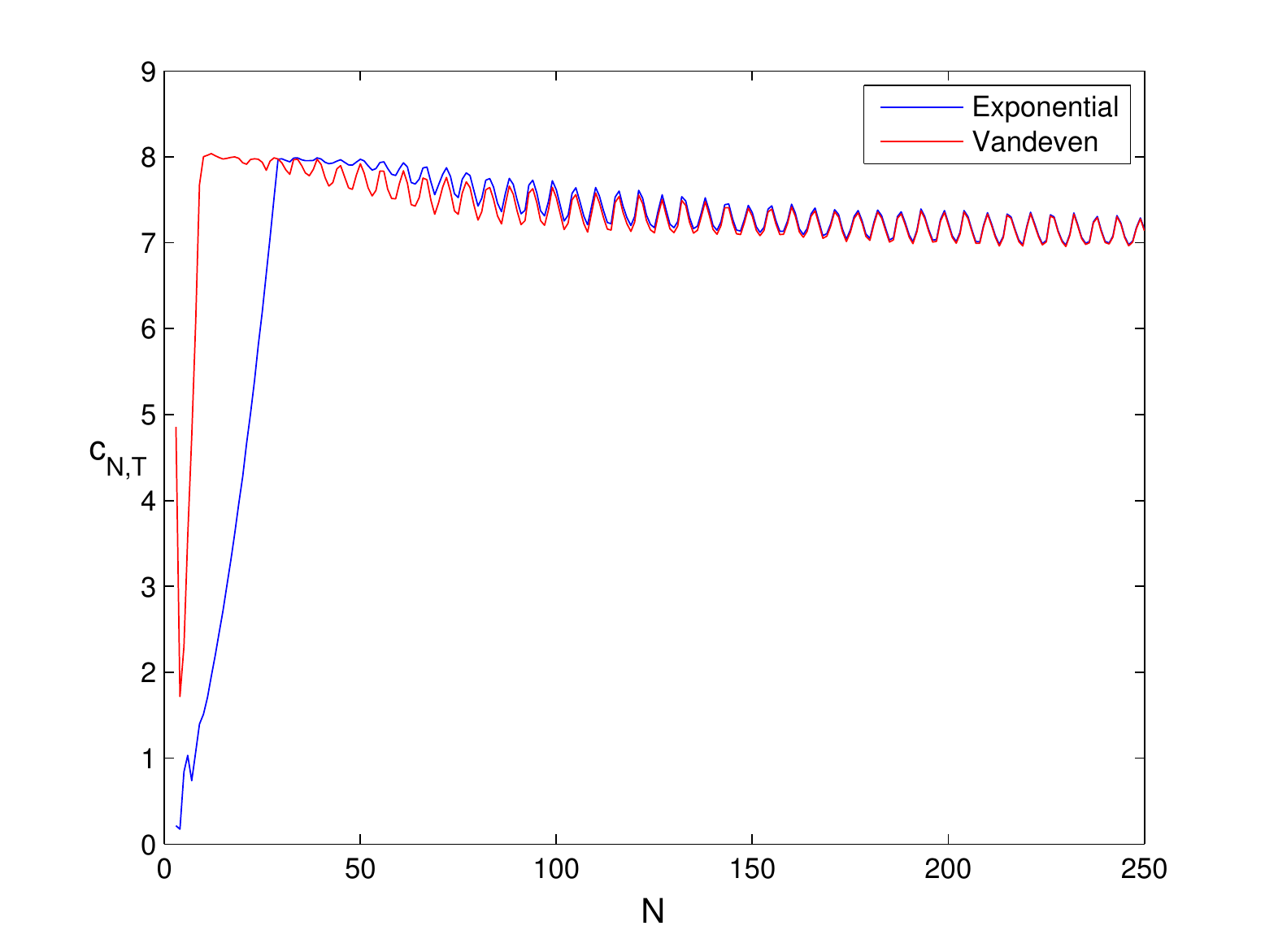}} 
 \mbox{\includegraphics[height=5cm]{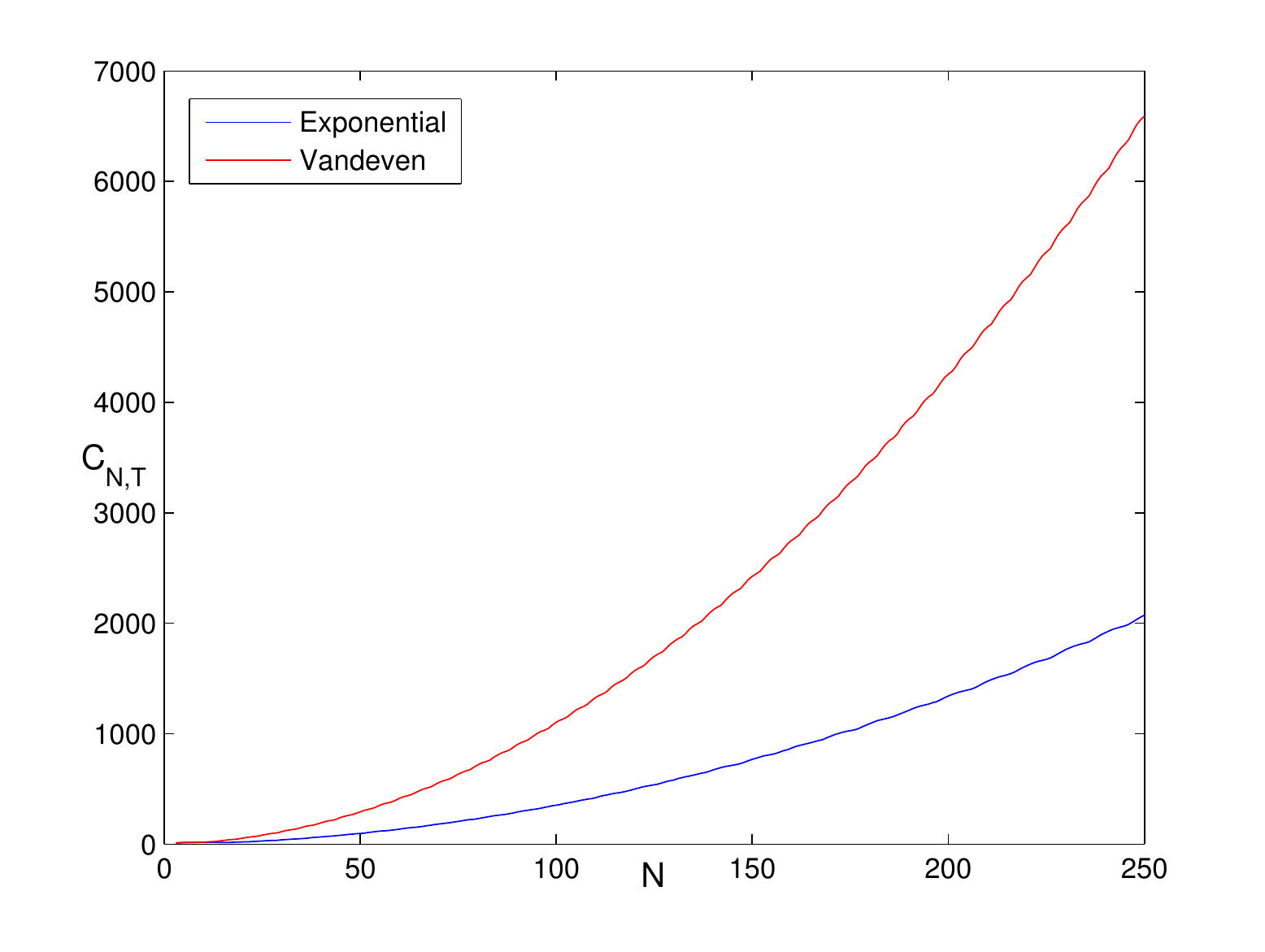}} \\
\end{tabular}
\end{center}
\caption{\footnotesize{Values of $c_{N,T}$ (left column) and $C_{N,T}$ (right column) for different filters and $T=8$. The exponential and Vandeven filters (bottom) are both of order $p=4$}}\label{const_filtre}
\end{figure}

Figure \ref{const_filtre} illustrates the behaviour of $c_{N,T}$ and $C_{N,T}$ of (\ref{obsfiltre})-(\ref{tracefiltre}) for different filters

\begin{figure}[!th]
\title{}
\begin{center}
\begin{tabular}{c}
 \mbox{\includegraphics[height=4cm]{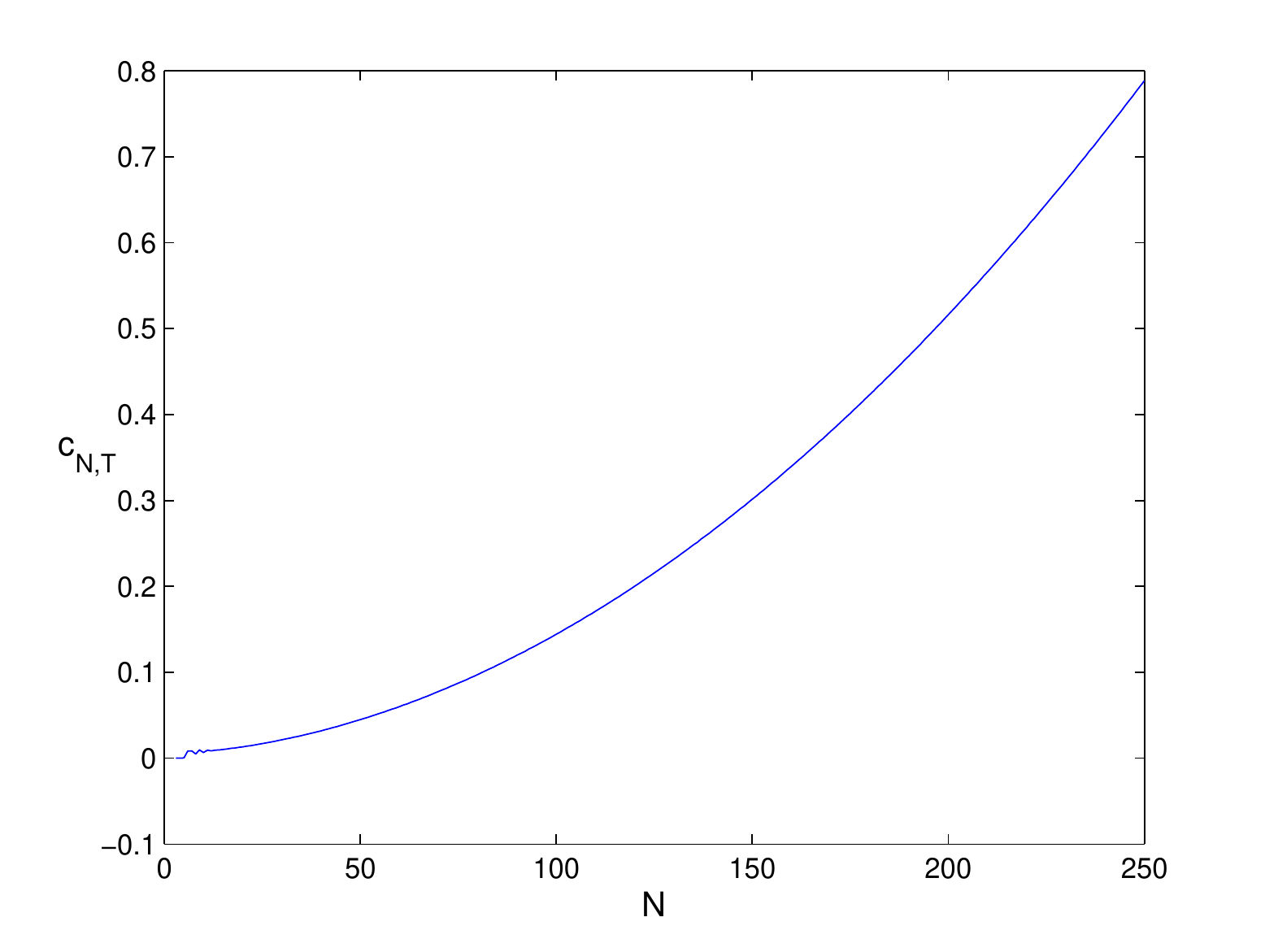}} 
 \mbox{\includegraphics[height=4cm]{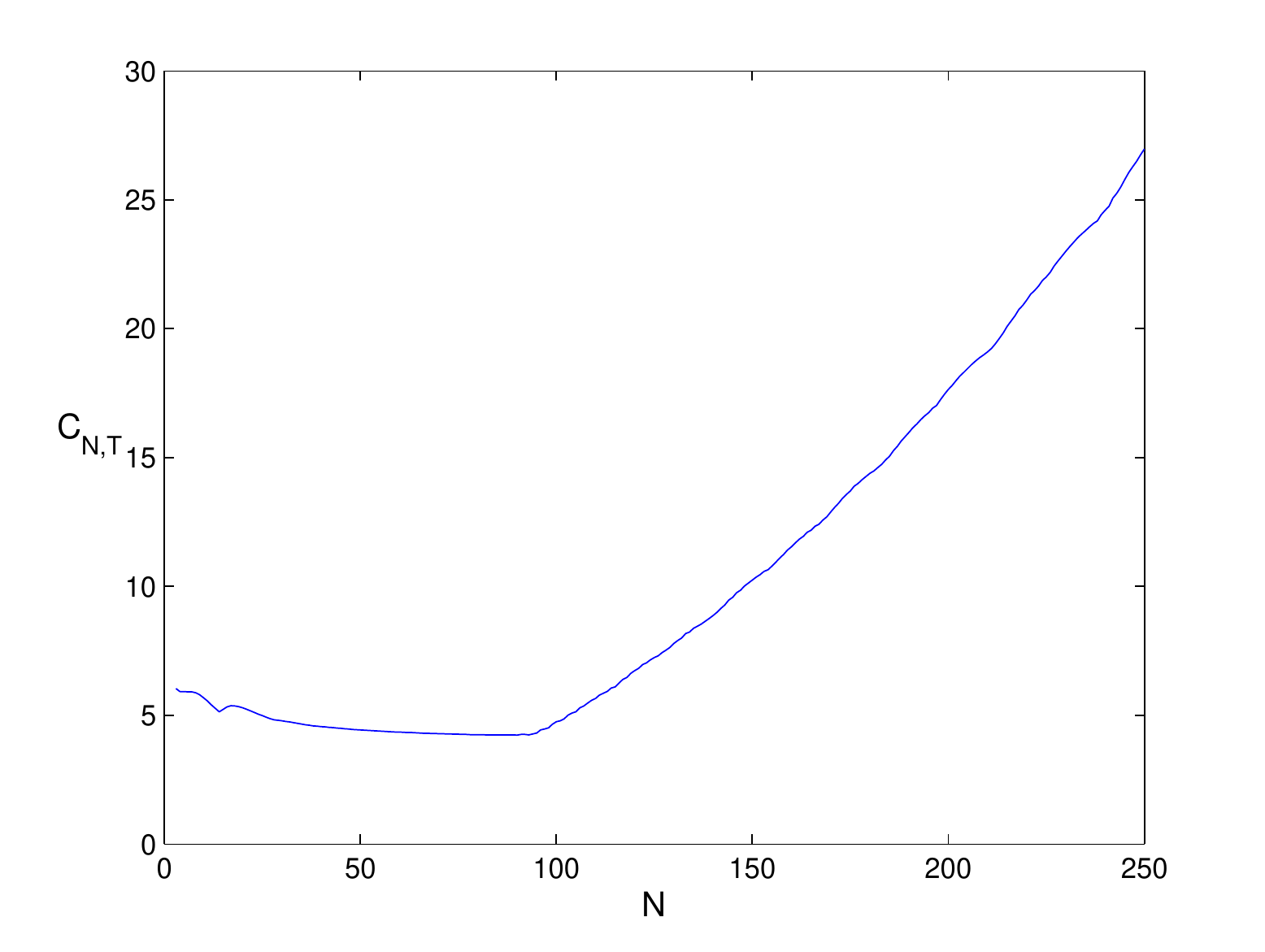}} \\
\end{tabular}
\end{center}
\caption{\footnotesize{Values of $c_{N,T}$ (left) and $C_{N,T}$ (right) for the exponential filter with $p=2$ and $T=8$}}\label{const_exp}
\end{figure}

According to these numerical results, there exists a uniform lower bound  on $c_{N,T}$ for the spectral filters considered here. However, we did not find any that could give an upper bound on $C_{N,T}$. Figure \ref{const_exp} shows our best result in this regard. A seemingly uniform bound on $C_{N,T}$ is obtained with the exponential filter but only up to a threshold value of $N$ (around $100$). Moreover, in the meantime, the lower bound for $c_{N,T}$ is dangerously low.

\section{Mixed Legendre-Galerkin Method}

As an alternative to the classical Legendre-Galerkin semi-discrete approximation of the homogeneous wave equation (\ref{obs}), we consider a mixed formulation, where $u$ and $u_t$ are the unknowns and are approximated directly (instead of $u$ alone in the classical formulation).

Let us denote $z=u_t$. Then  a mixed formulation is to find  $(u(t),z(t))\in H^1_0\times L^2$ satisfying 
\begin{align}
 \int_{-1}^1 u_t(x,t)\widetilde{v}(x)\, \textrm{dx}&=\int_{-1}^1 z(x,t)\widetilde{v}(x)\, \textrm{dx} & \forall &\widetilde{v}\in L^2  \label{Exactmixte1} \\
 \dfrac{d}{dt}<z(t),\hat{v}>_{H^{-1},H^1_0}&=-\int_{-1}^1 u_x(x,t)\hat{v}_x(x)\, \textrm{dx} & \forall  & \hat{v}\in H^1_0, \label{Exactmixte2}
\end{align}
for $0\leq t\leq T$ and initial conditions
\begin{eqnarray}
 u(0)=u_0\in H^1_0, \quad u_t(0)=u_1\in L^2.  \label{MixCI}
\end{eqnarray}

The mixed Legendre-Galerkin method considered here consists in seeking an approximation $(u^N(t),z^N(t))\in \PP^N_0\times \PP^{N-2}$, $0\leq t\leq T$, where \mbox{$\PP^M=\textrm{span}\{L_0(x),\ldots,L_M(x)\}$}, $M\geq 0$, satisfying 
\begin{align}
 \int_{-1}^1 u^N_t(x,t)\widetilde{v}^N(x)\, \textrm{dx}&=\int_{-1}^1 z^N(x,t)\widetilde{v}^N(x)\, \textrm{dx}  &\forall& \widetilde{v}^N\in \PP^{N-2}  \label{mixte1} \\
 \dfrac{d}{dt}<z^N(t),\hat{v}^N>_{H^{-1},H^1_0}&=-\int_{-1}^1 u^N_x(x,t)\hat{v}^N_x(x)\, \textrm{dx}  &\forall&  \hat{v}^N\in \PP^N_0, \label{mixte2}
\end{align}
and 
\begin{eqnarray}
 u^N(0)=u_0^N\in \PP^N_0 , \quad z(0)=u_1^N\in \PP^{N-2},  \label{MixAppCI}
\end{eqnarray}
where $u_0^N$ and $u_1^N$ are suitable approximations of $u_0$ and $u_1$ respectively. 
Writing down
\begin{eqnarray*}
u^N(x,t)&=&\sum_{k=1}^{N-1}a_{N,k}(t)\widetilde{L_k}(x) \\ 
z^N(x,t)&=&\sum_{k=0}^{N-2} b_{N,k}(t)L_k(x),
\end{eqnarray*}
and by letting $\textbf{a}_N(t)=:(a_{N,1}(t),\ldots,a_{N,N-1}(t))^t$ and $\textbf{b}_N(t)=:(b_{N,0}(t),\ldots,b_{N,N-2}(t))^t$, (\ref{mixte1})-(\ref{mixte2}) can be rewritten in the matrix form 
\begin{equation*}
\left(\begin{array}{c}
       \textbf{a}_{N}(t)' \\
       \textbf{b}_N(t)' 
      \end{array}
\right) 
=
\left(\begin{array}{cc}
      D_{N} & 0 \\
      0 & D_{N}^t \\
      \end{array}
\right)^{-1} 
\left(\begin{array}{cc}
      0 & M_N \\
      -K_N & 0 \\
      \end{array}
\right) 
\left(\begin{array}{c}
       \textbf{a}_N(t) \\
       \textbf{b}_N(t) 
      \end{array}
\right)
\end{equation*}
where $K_N\, \, M_N,\, D_{N}\in M_{N-1,N-1}(\RR)$ are given by  
\begin{equation*}
\begin{array}{rcl}
K_N(i,j) &=& \left\{ \begin{array}{ll}
        1 & \textrm{ if i=j} \\
        0 & \textrm{ otherwise}
       \end{array}\right. \\
M_N(i,j) &=& \left\{ \begin{array}{ll}
        \dfrac{2}{2i-1} & \textrm{ if i=j} \\
        0 & \textrm{ otherwise}
       \end{array}\right. \\
D_{N}(i,j) &=& \left\{ \begin{array}{ll}
       \dfrac{\sqrt{2}}{\sqrt{2i+1}(2i-1)} & \textrm{ if i=j} \\
       \dfrac{-\sqrt{2}}{\sqrt{2j+1}(2j+3)} & \textrm{ if i=j+2} \\
       0 & \textrm{ otherwise}.
      \end{array}\right.
\end{array}
\end{equation*}
Note that 
\begin{equation*}
 \left(\begin{array}{cc}
      D_{N} & 0 \\
      0 & D_{N}^t \\
      \end{array}
\right)
\end{equation*}
is nonsingular.

%
%
%

\begin{figure}[!th]
\begin{center}
\begin{tabular}{c}
 \mbox{\includegraphics[height=4cm]{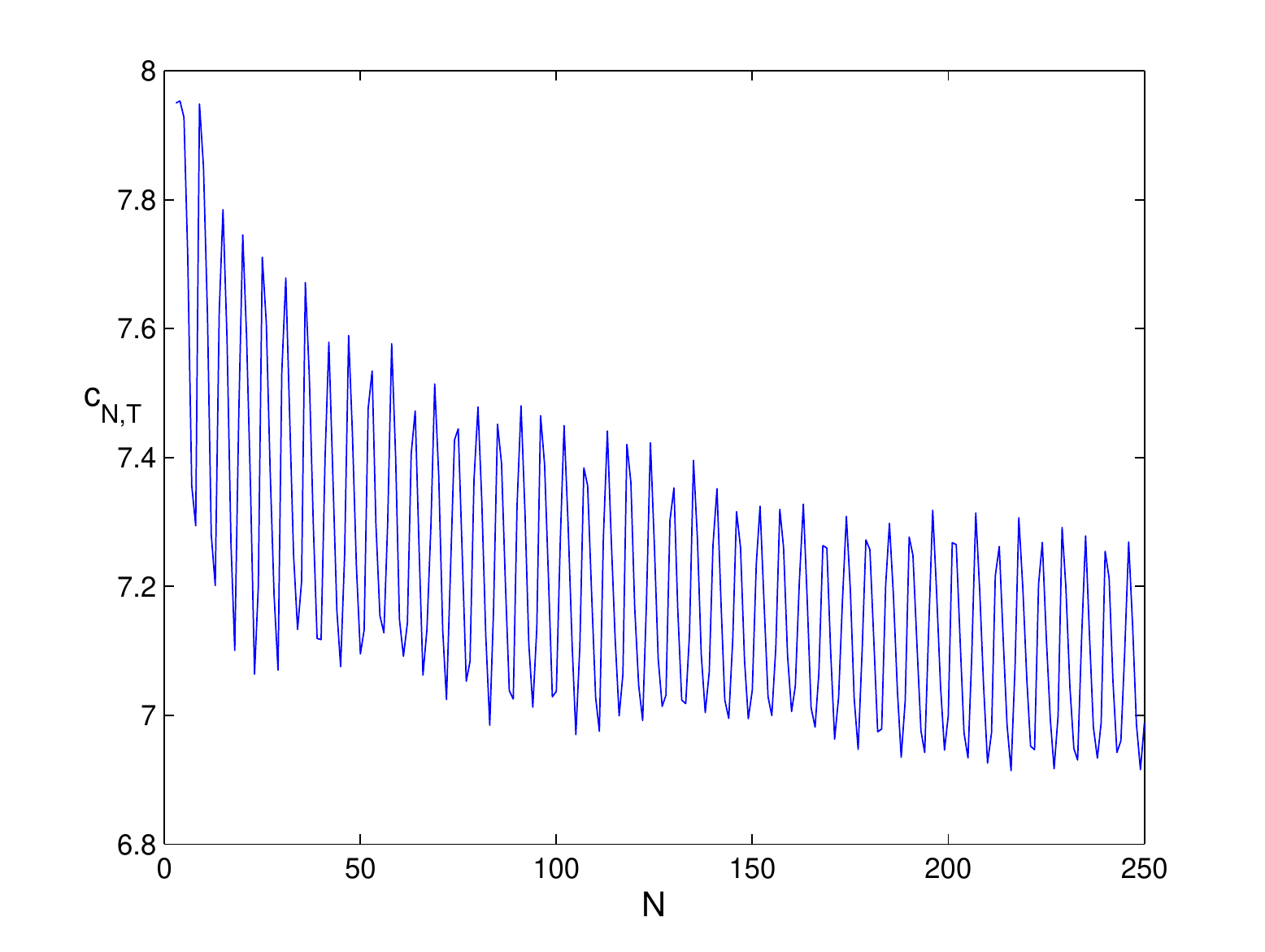}} 
 \mbox{\includegraphics[height=4cm]{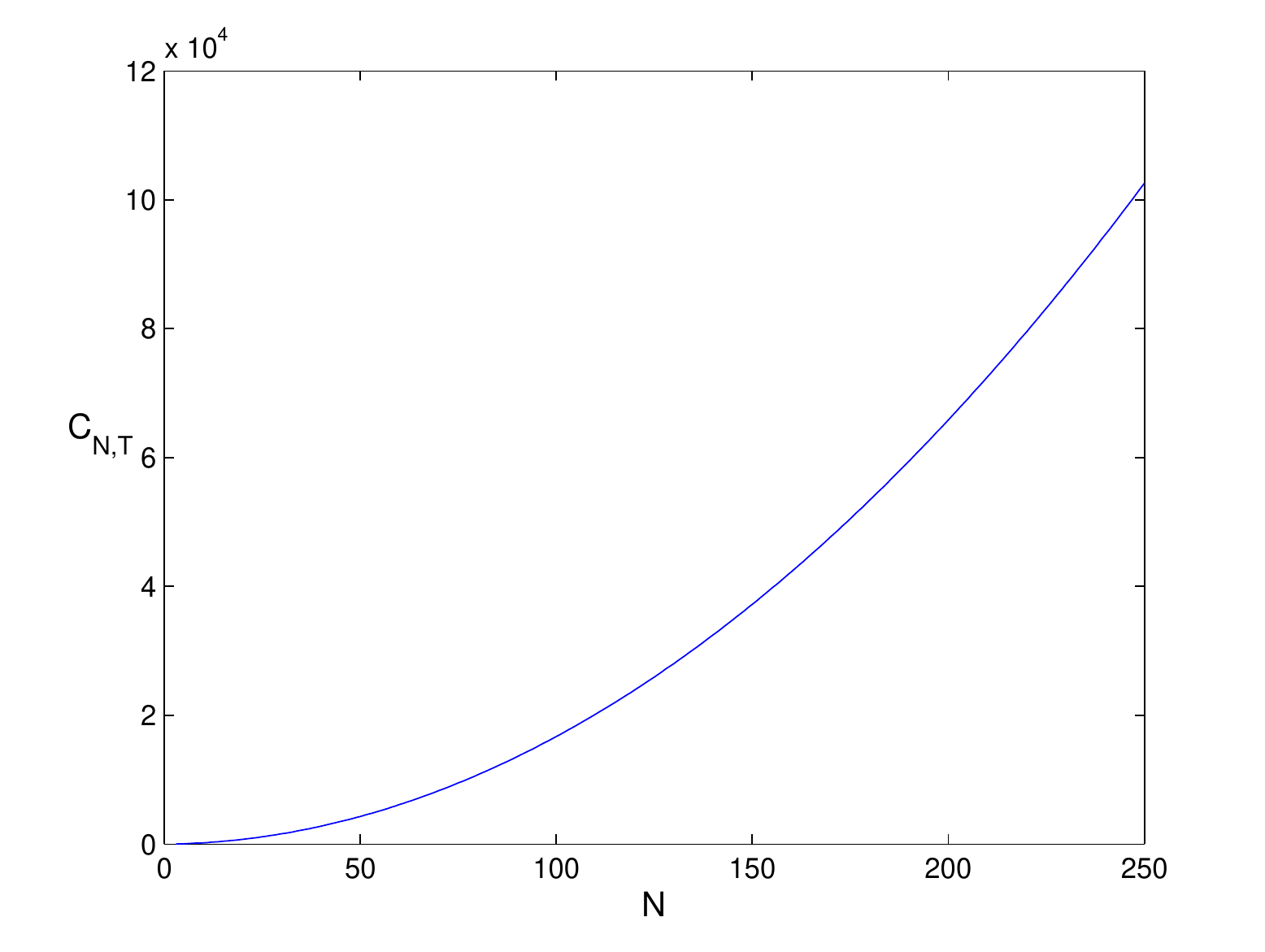}} 
\end{tabular}
\end{center}
\caption{\footnotesize{Values of $c_{N,T}$ (left) and $C_{N,T}$ (right) associated to the mixed Legendre Galerkin formulation (\ref{mixte1})-(\ref{mixte2}) with $T=8$}}\label{constmixte}
\end{figure}

Here again, the energy of the solutions of (\ref{mixte1})-(\ref{mixte2}) is conserved along time
\begin{equation*}
E(u^N(t),z^N(t))=E(u^N(0),z^N(0)),
\end{equation*}
where
\begin{equation*}
E(u^N(t),z^N(t))=\dfrac{1}{2}\int_{-1}^1 \left|z^N(x,t)\right|^2 + \left|u^N(x,t)\right|^2 \textrm{dx}.
\end{equation*}
Let us now define the observability and continuity constants $c_{N,T}$ and $C_{N,T}$, respectively, by   

\begin{equation*}
c_{N,T} = \inf_{(u_0^N,u_1^N)\in {\PP_N^0\times\PP_{N-2}}} \left(\int_0^T \left| u^N_x(1,t)\right|^2 dt \right) \;/ \; E(u^N(0),z^N(0)),
\end{equation*}

\begin{equation*}
C_{N,T} = \sup_{(u_0^N,u_1^N)\in {\PP_N^0\times\PP_{N-2}}} \left(\int_0^T \left| u^N_x(1,t)\right|^2 dt \right) \;/ \; E(u^N(0),z^N(0)),
\end{equation*}
Their behaviour are shown  in Figure \ref{constmixte}. We see that we recover a positive uniform lower bound for $c_{N,T}$, but no uniform upper bound for $C_{N,T}$. 

\section{Nitsche's method}

Nitsche's method \cite{Nitsche} was first introduced to weakly impose boundary conditions of Dirichlet type. Subsequently, the method was used to (weakly) impose the inter-element continuity in the discontinuous Galerkin finite element method (DGFEM) \cite{ARN}. The DGFEM combined with Nitsche's method have been applied to the second order wave equation \cite{GROTE} and even studied for the numerical boundary observability of waves \cite{MaricaZuazua}. Note that they do not provide a uniform numerical boundary observability, unless Fourier filtering or a multi-grid strategy is applied. 

Nitsche's method has similarities with the penalty method but unlike the latter it is consistent and does have its ill-conditioning problems at the discrete level. In fact, it has even more similarities with a stabilized Lagrange multiplier method even if it does not involve any multiplier \cite{hansbo}.  

Here we consider Nitsche's method for the weak imposition of the Dirichlet boundary condition at $x=1$, which is to say at the end point where the solutions are observed or controlled. We also performed numerical experiments with Nitsche's method used at both endpoints, but the results are qualitatively the same. 

There is two versions of Nitsche's method. A symmetrical one and an unsymmetrical one. Both depend on a parameter which, for the stability of the discrete formulation, has to be chosen carefully, depending on the approximation space, among other things.  We present the two versions and the corresponding numerical results successively, beginning with the symmetrical case which, unlike the unsymmetrical case, has the property of still leading to a conservative finite dimensional dynamical system. 

\subsection{Symmetric case}

Let  $H^1_L(-1,1)=\{ \psi\in H^1(-1,1) \, | \, \psi(-1)=0 \}$ and let $\PP^N_L$ be the set of polynomials of degree $N$ vanishing at $x=-1$. We have, for instance, 
\begin{equation*}
	\PP^N_L=\textrm{span}\{\widehat{L_1}(x),\ldots,\widehat{L_N}(x)\}\subset H^1_L(-1,1), 
\end{equation*}
where $\widehat{L_k}(x)=L_k(x)-(-1)^{k}$, $k\geq 1$.
Let us define the following bilinear form over $\PP^N_L \times \PP^N_L$:
\begin{eqnarray*}
A_N(u^N,\psi^N)&=&\int_{-1}^1 u^N_{x}(x,t)\psi_x^N(x) \,\textrm{dx} - u^N_x(1,t)\psi^N(1)  \\
&& -u^N(1,t)\psi^N_x(1)+\gamma N^2u^N(1,t)\psi^N(1). 
\end{eqnarray*}
Note that this bilinear form is symmetric. The approximation method using the symmetric variant of Nitsche's method then reduces to seeking an approximation $u^N(.,t)\in\PP^N_L$, $t>0$, for instance in the form 
\begin{equation*}
 u^N(x,t)=\sum_{k=1}^{N} a_{N,k}(t)\widehat{L_k}(x),
\end{equation*}
satisfying 
\begin{equation}\label{Nitsche}
\int_{-1}^1 u^N_{tt}(x,t)\psi^N(x) \, \textrm{dx} + A_N(u^N,\psi^N)=0, \, \forall v^N\in \PP^N_L.
\end{equation}
The linear system (\ref{Nitsche}) can be written in the matrix form  
\begin{eqnarray*}
M_N\textbf{a}_N''(t)+K_{N}\textbf{a}_N(t)=0 
\end{eqnarray*}
where $\textbf{a}_N(t)=(a_{N,1}(t),\ldots,a_{N,N}(t))^t$ and  $(M_N,K_{N})\in M_{N\times N}(\RR)^2$ are given by
\begin{equation*}
\begin{array}{rcl}
M_N(i,j)&=&\left\{\begin{array}{l l}\frac{4i+4}{2i+1} & \textrm{ if i=j}\\
         -2 & \textrm{ if } i+j \textrm{ odd}\\
	 2 & \textrm{ if } i+j \textrm{ even and } i\neq j
         \end{array}\right.\\ 
K_{N}&=&P_{N}-Q_{N}-Q_{N}^t+\gamma N^2 R_{N},\\
P_{N}(i,j)&=&\left\{\begin{array}{l l} r^2+r & \textrm{ if } i+j \textrm{ even}\\
	 0 & \textrm{ otherwise}
         \end{array}\right.\\
Q_{N}(i,j)&=&\left\{\begin{array}{l l} j^2+j & \textrm{ if } i \textrm{ odd}\\
	 0 & \textrm{ otherwise}
         \end{array}\right.\\
R_{N}(i,j)&=&\left\{\begin{array}{l l} 4 & \textrm{ if } i,j \textrm{ are odds}\\
	 0 & \textrm{ otherwise}
         \end{array}\right.\\
\end{array}
\end{equation*}

Analogously to what happens with convergence proofs of approximations based on the discontinuous Galerkin method using Nitsche's method (see \cite{GROTE}), we need the bilinear form $A_N$ to be coercive (with respect to a $H^1_L(-1,1)$-norm) in $\PP^N_L \times \PP^N_L$.  This depends on the value of $\gamma$, and we show next that coercivity is ensured with $\gamma>1/2$.  In order to prove the coercivity under this condition, let us recall the following inverse inequality (see \cite{WAR}) : 
\begin{equation*}
|p^N(1)|^2\leq \dfrac{(N+1)^2}{2} |p^N(x)|^2_{H^1_0(-1,1)}, \ \forall p^N(x)\in\PP^N(-1,1)
\end{equation*}
Then, the coercivity easily follows, provided $\gamma>1/2$:
\begin{eqnarray*}
A_N(\psi^N,\psi^N)&=& |\psi^N|_{H^1_0}^2 -2\psi^N_x(1)\psi^N(1)+\gamma N^2\left|\psi^N(1)\right|^2 \\
            &\geq& |\psi^N|_{H^1_0}^2 -2\left|N^{-1}\psi^N_x(1)\right|\left|N\psi^N(1)\right|+\gamma N^2\left|\psi^N(1)\right|^2 \\
	    &\geq& |\psi^N|_{H^1_0}^2 -\frac{1}{\epsilon}\left|N^{-1}\psi^N_x(1)\right|^2 - \epsilon \left|N\psi^N(1)\right|^2+\gamma N^2\left|\psi^N(1)\right|^2, \quad \forall \epsilon > 0 \\
	    &\geq& (1-\frac{1}{2\epsilon})|\psi^N|_{H^1_0}^2 + (\gamma-\epsilon)\left|N \psi^N(1)\right|^2 , \quad \forall \epsilon > 0\\
	    &\geq& C\|\psi^N\|_{1,N}^2,
\end{eqnarray*}
for $C>0$ if $\gamma>\epsilon>1/2$ and  where $\|.\|_{1,N}$ denotes the following norm on $H^1_L(-1,1)$
\begin{equation*}
\|\psi\|^2_{1,N}:=\left|\psi\right|^2_{H^1_0(-1,1)}+N^2\left|\psi(1)\right|^2.
\end{equation*}
Note that the coercivity constant $C$ does not depend on $N$ and that $A_N$ is coercive, also uniformly with respect to $N$, if instead of $\|.\|^2_{1,N}$ we take the equivalent norm $|.|_{H^1_0}$.


Taking $\psi^N=u_t^N(x,t)$ in (\ref{Nitsche}), we obtain that the energy of the numerical solutions of (\ref{Nitsche}), defined by   
\begin{eqnarray*}
\lefteqn{E_{N,\gamma}(u^N(t))=\dfrac{1}{2}\int_{-1}^1 \left|u^N_t(x,t)\right|^2 + \left|u^N_x(x,t)\right|^2 \textrm{dx}}\nonumber  \\
& &\qquad + \left(\dfrac{\gamma}{2} N^2u^N(1,t)-u^N_{x}(1,t)\right)u^N(1,t),
\end{eqnarray*}
is preserved along time : $E_{N,\gamma}(u^N(t))=E_{N,\gamma}(u^N(0))$. Note that we have proved that if $\gamma>1/2$ then $E^{1/2}_{N,\gamma}$ is a norm in the energy space $H^1_L(-1,1)\times L^2(-1,1)$ for solutions of the wave equation and solutions of (\ref{Nitsche}) as well. 

Regarding our observability problem, we define the following discrete observability and continuity constants $c_{N,T}$ and $C_{N,T}$ respectively:
\begin{equation}\label{obsnitsche}
c_{N,T} = \inf_{(u_0^N,u_1^N)\in {(\PP^N_L)}^2} \left(\int_0^T \left| u^N_x(1,t)-\gamma N^2 u^N(1,t)\right|^2 dt \right) \;/ \; E_{N,\gamma}(u^N(0)),
\end{equation}

\begin{equation}\label{contnitsche}
C_{N,T} = \sup_{(u_0^N,u_1^N)\in {(\PP^N_L)}^2} \left(\int_0^T \left| u^N_x(1,t)-\gamma N^2 u^N(1,t)\right|^2 dt \right) \;/ \; E_{N,\gamma}(u^N(0)).
\end{equation}


Some remarks concerning these choices are in order since the observation is not $u_x^N(1,t)$ any more, but $u^N_x(1,t)-\gamma N^2 u^N(1,t)$ instead. 

First, if we use Nitsche's method and the resulting discrete formulation (\ref{Nitsche}) only to have the approximation $u^N$ involved in the functional $J_N$ (defined in (\ref{mindisc})) to be minimized then it is natural to observe  $u_x^N(1,t)$. Unfortunately, our numerical results show that if we drop the term  $\gamma N^2 u^N(1,t)$ then  $c_{N,T}$ is not uniformly bounded from below by a positive constant. 

Alternatively, if we also use Nitsche's method to approximate the controlled wave equation, that is to append the Dirichlet boundary condition $y(1,t)=v(t)$, then we obtain that the associated observability and continuity constants are defined by (\ref{obsnitsche}) and (\ref{contnitsche}), and that the controllability problem reduces to minimize the functional 
\begin{eqnarray*}
J^N_{\gamma}(u^N_0,u^N_1)&=&\dfrac{1}{2}\int_0^T \left| -u^N_x(1,t)+\gamma N^2u^N(1,t) \right|^2 \textrm{dt}
\\ && \quad -<(y_1,-y_0),(u^N_0,u^N_1)>_{H^{-1}\times L^2,H^1_0\times L^2}.
\end{eqnarray*}
over $(\PP^N_L)^2$. 

Indeed, let  $y^N(.,t)\in\PP^N_L$ be the numerical approximation of $y(.,t)$ satisfying  
\begin{eqnarray}\label{Nitschey}
\int_{-1}^1 y^N_{tt}(x,t)\psi^N(x) \, \textrm{dx} + \int_{-1}^1 y^N_{x}(x,t)\psi_x^N(x) \, \textrm{dx} - y^N_x(1,t)\psi^N(1)-y^N(1,t)\psi^N_x(1) \nonumber \\
 +\gamma N^2y^N(1,t)\psi^N(1)=v^N(t)\left(-\psi^N_x(1)+\gamma N^2\psi^N(1)\right), \quad \forall \psi^N\in \PP^N_L. \qquad 
\end{eqnarray}
Let's also consider, for the moment, that $v^N(t)$ is a control, that is $y^N(x,T)=y^N_t(x,T)=0$. Now let, in (\ref{Nitschey}), $\psi^N=u^N$, where $u^N(x,t)$ is the solution of (\ref{Nitsche}), and let's integrate this relation in time. By integrating by parts twice in time the first integral, we obtain
\begin{eqnarray*}
\int_0^T \int_{-1}^1 y^N(x,t)u^N_{tt}(x,t) \, \textrm{dxdt} + \int_0^T \int_{-1}^1 y^N_{x}(x,t)u_x^N(x,t) \,\textrm{dxdt} -\int_0^T  y^N_x(1,t)u^N(1,t) \, \textrm{dt} \\ 
-\int_0^T y^N(1,t)u^N_x(1,t)  \, \textrm{dt} +\gamma N^2 \int_0^T y^N(1,t)u^N(1,t) \, \textrm{dt} \qquad \qquad \qquad \quad  \\
 =\int_0^T v^N(t)\left(-u^N_x(1,t)+\gamma N^2u^N(1,t)\right)  \, \textrm{dt} + \int_{-1}^1 u^N_0(x)y_1^N(x)-u^N_1(x)y_0^N(x) \, \textrm{dx},  \qquad 
\end{eqnarray*}
The left-hand side corresponds to the left-hand side of (\ref{Nitsche}) where $\psi^N=y^N$. Thus, we obtain
\begin{equation}
 \int_0^T v^N(t)\left(-u^N_x(1,t)+\gamma N^2u^N(1,t)\right)  \, \textrm{dt} = \int_{-1}^1 u^N_1(x)y_0^N(x)- u^N_0(x)y_1^N(x)\, \textrm{dx}.
\end{equation}
Using the same arguments as in the introduction for the continuous wave equation, we obtain that the approximate control of minimal $L^2(0,T)$-norm is then $v^N(t)=u^N_x(1,t)-\gamma N^2 u^N(1,t)$ where $u^N$ is the solution of \eqref{Nitsche} with the initial data $(u_0^N,u_1^N)$ minimizing $J^N_{\gamma}$ over $(\PP^N_L)^2$. Let us recall that a uniform positive lower bound $c_{N,T}\geq c_T>0$ ensures that $J^N_{\gamma}$ is uniformly coercive. 
%

\begin{figure}[!th]
\begin{center}
\begin{tabular}{c}
 \mbox{\includegraphics[height=4cm]{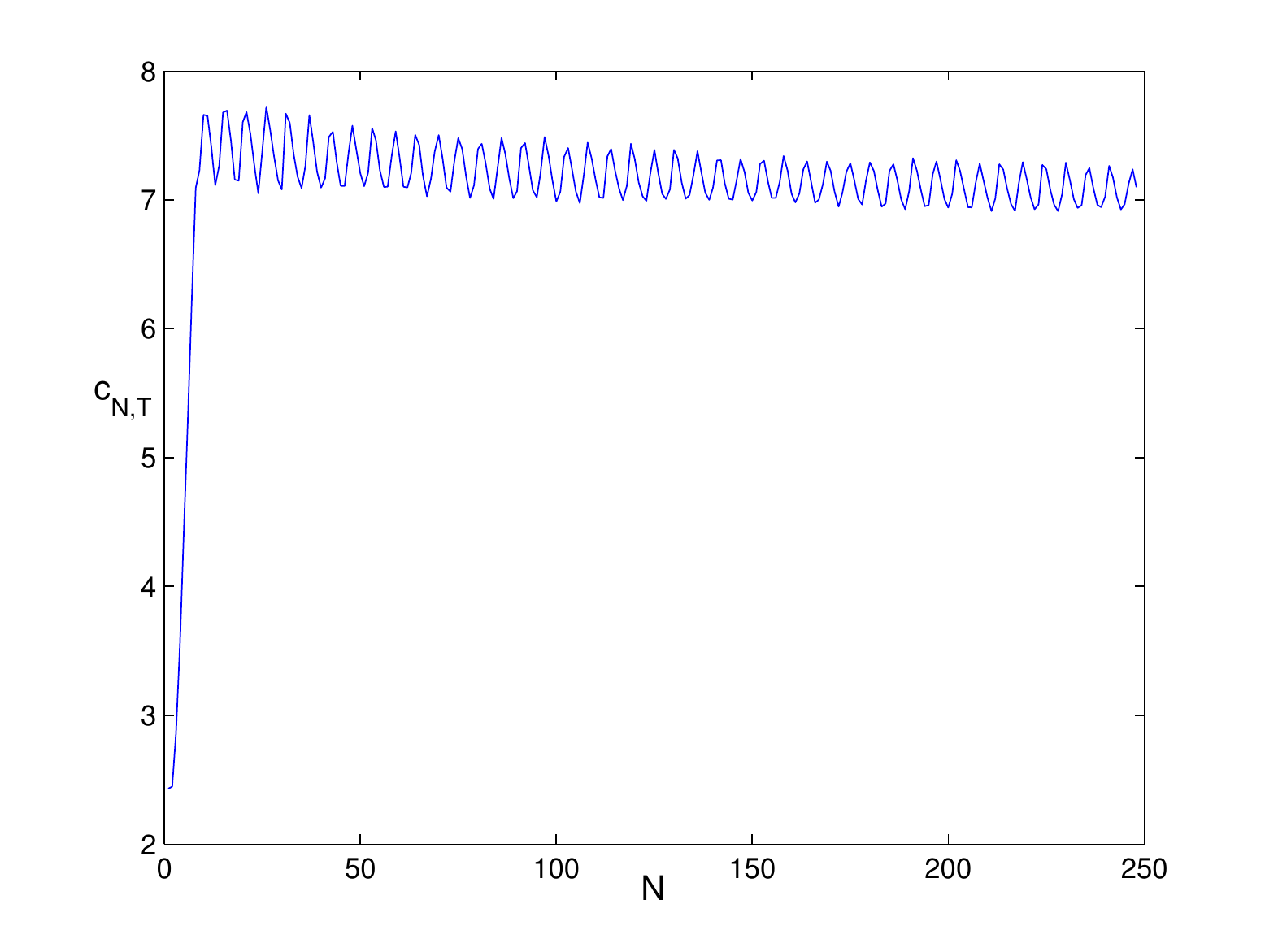}} 
 \mbox{\includegraphics[height=4cm]{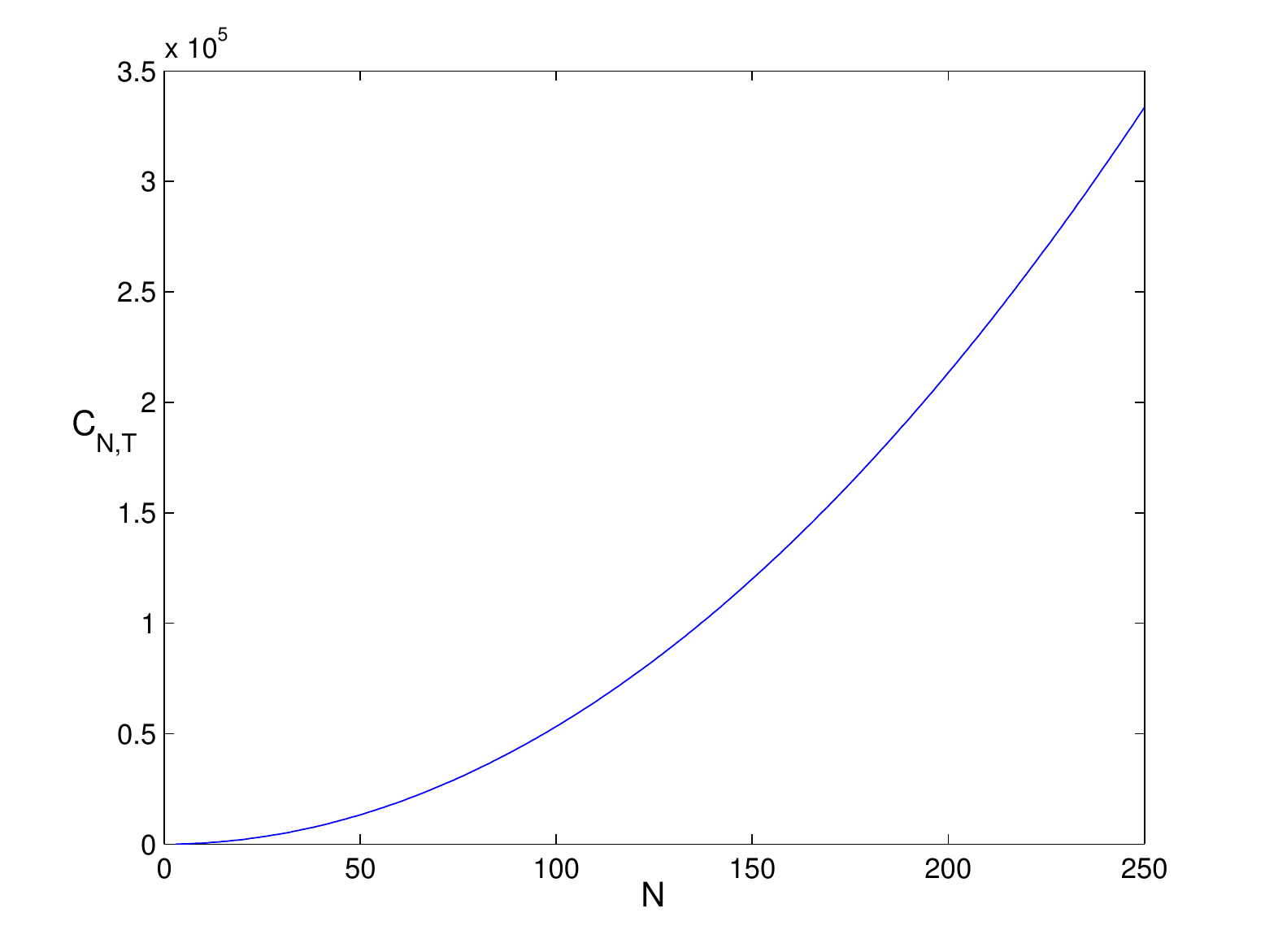}} 
\end{tabular}
\end{center}
\caption{\footnotesize{Values of $c_{N,T}$ (left) and $C_{N,T}$ (right) with the symmetric Nitsche's method, with $\gamma=0.8$ and $T=8$}}\label{constnitsche08}
\end{figure}

In Figure \ref{constnitsche08} we show the behaviour of $c_{N,T}$ and $C_{N,T}$ in (\ref{obsnitsche})-(\ref{contnitsche}) for the particular case $\gamma=0.8$. We obtain numerically an uniform lower bound for $c_{N,T}$ but no upper bound for $C_{N,T}$. 

%
%

\begin{figure}[!th]
\begin{center}
\begin{tabular}{c}
 \mbox{\includegraphics[height=4cm]{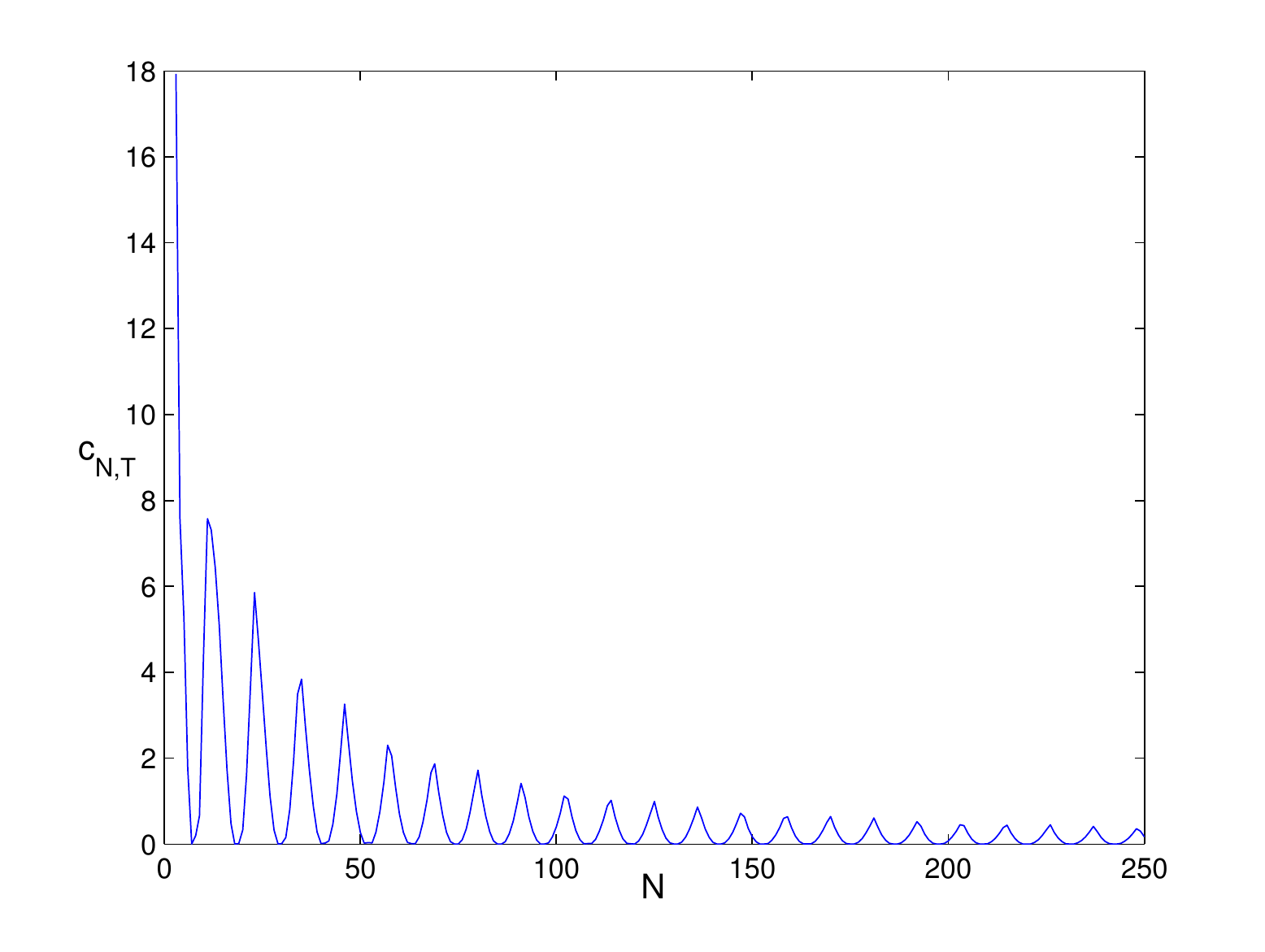}} 
 \mbox{\includegraphics[height=4cm]{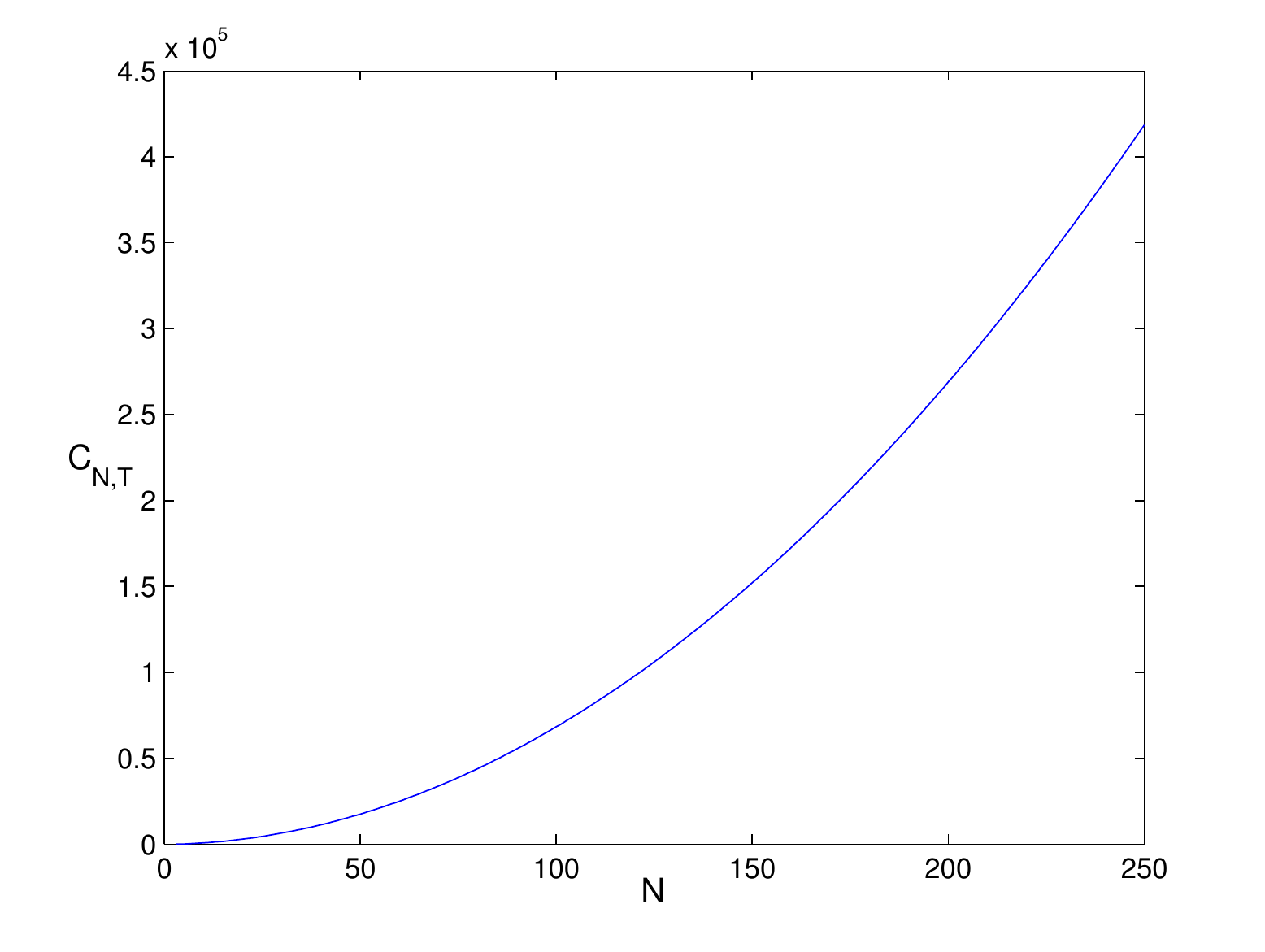}} 
\end{tabular}
\end{center}
\caption{\footnotesize{Values of $c_{N,T}$ (left) and $C_{N,T}$ (right) for the symmetric Nitsche's method, when the term $\gamma N^2u^N(1,t)$ is dropped in their definitions (\ref{obsnitsche}) and (\ref{contnitsche}),  with $\gamma=0.8$ and $T=8$}}\label{comparaisonNitschenonunif}
\end{figure}

In order to show the influence of the term $\gamma N^2u^N(1,t)$ in the observation and the resulting definition of $c_{N,T}$ and $C_{N,T}$, we show in Figure \ref{comparaisonNitschenonunif} that if we drop this term then, in the same situation than in the previous example ($T=8$, $\gamma=0.8$), we do not have uniform observability (positive lower bound on $c_{N,T}$) any more, nor we have uniform continuity of the trace (upper bound on $C_{N,T}$).

\subsection{Non-symmetric case}

Application of the non-symmetric Nitsche's method reduces to seek $u^N(.,t)\in \PP^N_L$ such that 
\begin{eqnarray}\label{Nitschenonsym}
\lefteqn{\int_{-1}^1 u^N_{tt}(x,t)\psi^N(x) \textrm{dx} + \int_{-1}^1 u^N_{x}(x,t)\psi_x^N(x) \textrm{dx} - u^N_x(1,t)\psi^N(1)} \nonumber \\
& & \qquad +u^N(1,t)\psi^N_x(1)+\gamma N^2u^N(1,t)\psi^N(1)=0, \forall \psi^N\in \PP^N_L.
\end{eqnarray}
where $\gamma$ is still a strictly positive parameter at our disposal. The interest of this non-symmetric formulation is that there is no condition on $\gamma$ to have a coercive bilinear form. Let us define the bilinear form in question by 
\begin{eqnarray*}
A^N(u^N,\psi^N)&:=&\int_{-1}^1 u^N_{x}(x,t)\psi_x^N(x)\, \textrm{dx} - u^N_x(1,t)\psi^N(1)+u^N(1,t)\psi^N_x(1)\\ &&\,+\ \gamma N^2u^N(1,t)\psi^N(1), 
\end{eqnarray*}  
one directly has that 
\begin{equation*}
A^N(\psi^N,\psi^N)=\|\psi^N\|^2_{1,\sqrt{\gamma}N}.
\end{equation*}

Let $\textbf{a}_N(t)=(a_{N,1}(t),\ldots,a_{N,N}(t))^t$. Then, (\ref{Nitschenonsym}) can be written in the matrix form \begin{eqnarray*}
M_N\textbf{a}_N''(t)+K_N\textbf{a}_N(t)=0 
\end{eqnarray*}
where  $\left(M_N,K_N\right)\in M_{N\times N}(\RR)^2$, $K_N=P_N-Q_N+Q_N^t+\gamma N^2 R_N$ and where $M_N$, $P_N$, $Q_N$ and $R_N$ were given previously.   

The loss of symmetry of the bilinear elliptic form makes us  consider the following definitions for the observability and continuity constants:

\begin{equation}\label{NitscheNonSymobs}
c_{N,T} = \inf_{(u_0^N,u_1^N)\in {(\PP^N_L)}^2} \left(\int_0^T \left| u^N_x(1,t)+\gamma N^2 u^N(1,t)\right|^2 dt \right) \;/ \; E(u^N(0)),
\end{equation}
\begin{equation}\label{NitscheNonSymtrace}
C_{N,T} = \sup_{(u_0^N,u_1^N)\in {(\PP^N_L)}^2} \left(\int_0^T \left| u^N_x(1,t)+\gamma N^2 u^N(1,t)\right|^2 dt \right) \;/ \; E(u^N(0)).
\end{equation}
where $E(.)$ stands for the energy of the continuous wave equation (\ref{defenergie}). Even if the energy of the solutions of (\ref{Nitschenonsym}) is not preserved along time, $E(.)$ still is a norm in the energy space $H^1_L(-1,1)\times L^2(-1,1)$. Note that the observation $u^N_x(1,t)+\gamma N^2 u^N(1,t)$ was obtained by the exact same way than in the previous subsection, by considering the approximation of the controlled wave equation using the same non-symmetrical Nitsche's method.  

\begin{figure}[!th]
\begin{center}
\begin{tabular}{c}
 \mbox{\includegraphics[height=4cm]{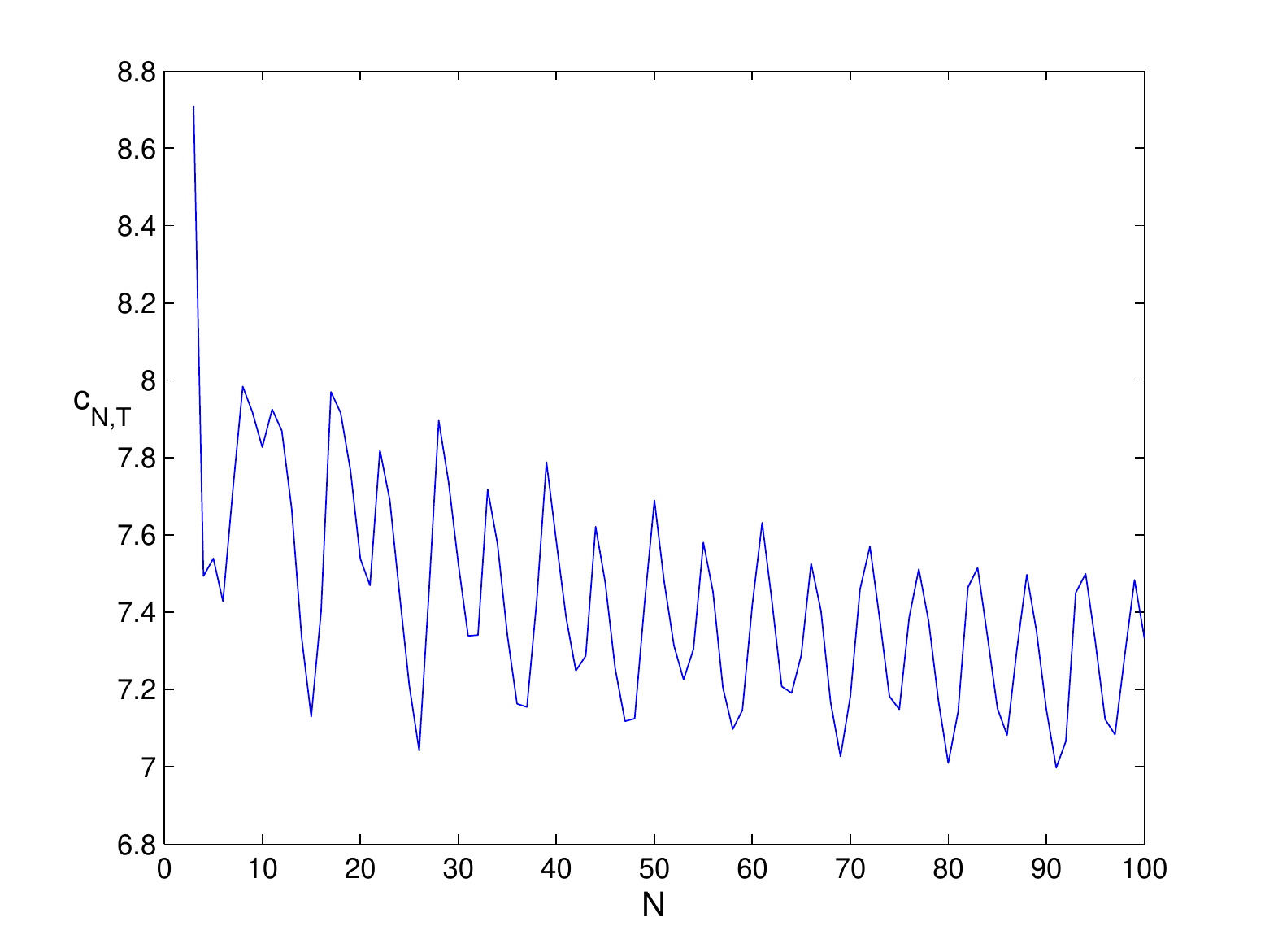}} 
 \mbox{\includegraphics[height=4cm]{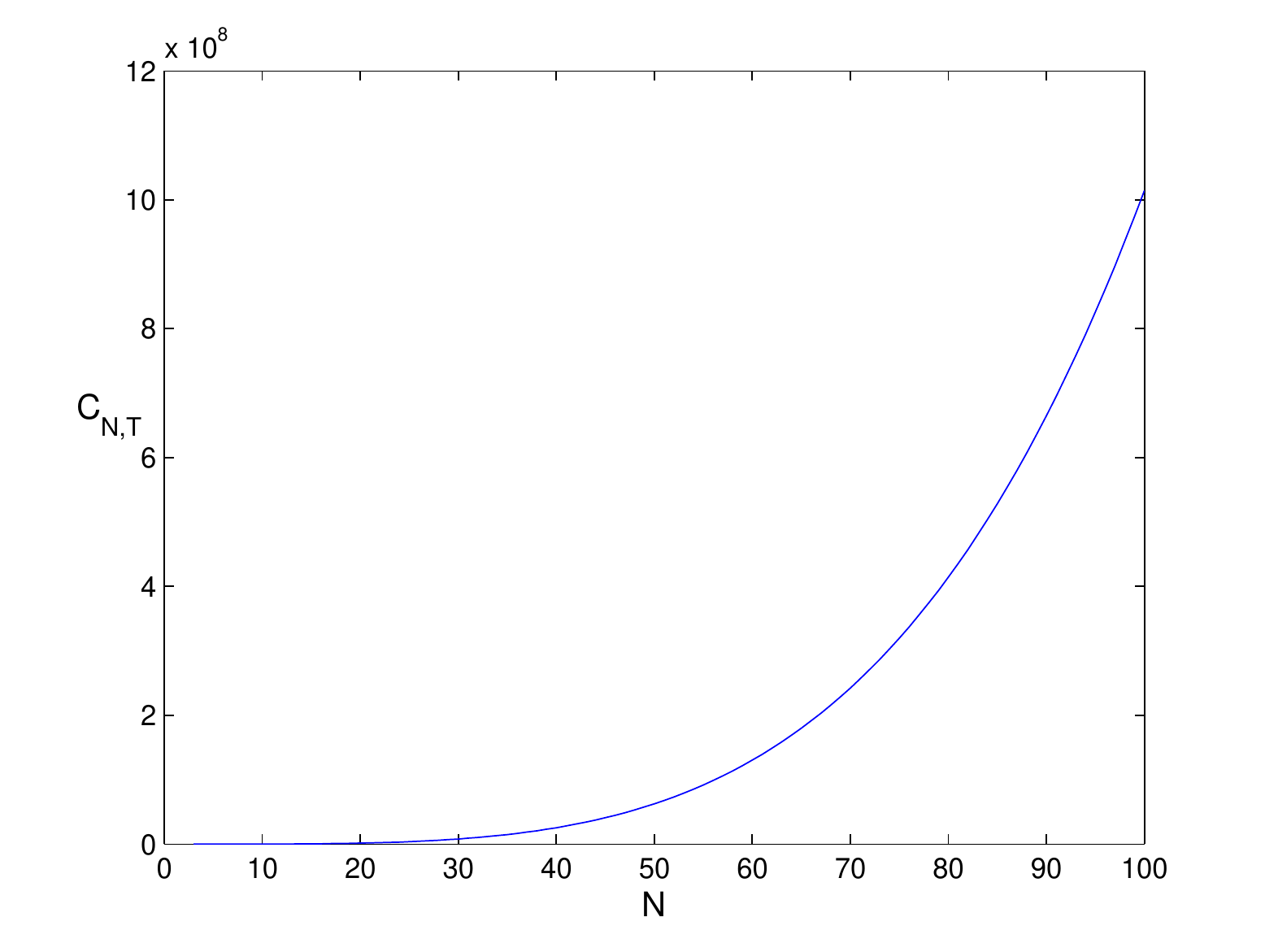}} 
\end{tabular}
\end{center}
\caption{\footnotesize{Values of $c_{N,T}$ (left) and $C_{N,T}$ (right) defined by (\ref{NitscheNonSymobs}) and (\ref{NitscheNonSymtrace}), for the non-symmetric Nitsche's method with $\gamma=1$ and $T=8$}}\label{constns1}
\end{figure}

In Figure \ref{constns1} we show the behaviour of $c_{N,T}$ and $C_{N,T}$ for $\gamma=1$ and $T=8$. Again, we observe a uniform observability but no upper bound for $C_{N,T}$ .
%

\begin{figure}[!th]
\begin{center}
\begin{tabular}{c}
 \mbox{\includegraphics[height=4cm]{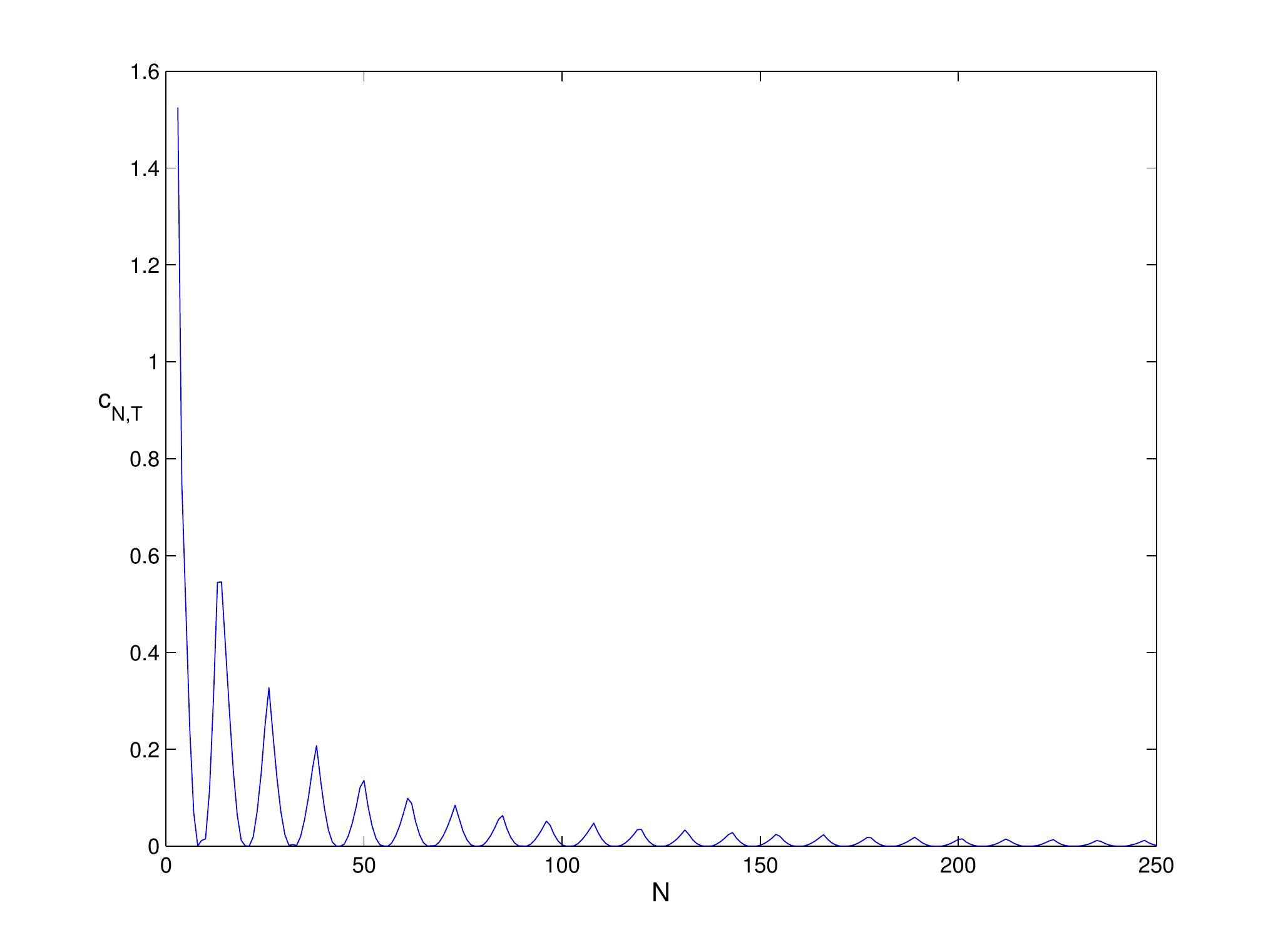}} 
 \mbox{\includegraphics[height=4cm]{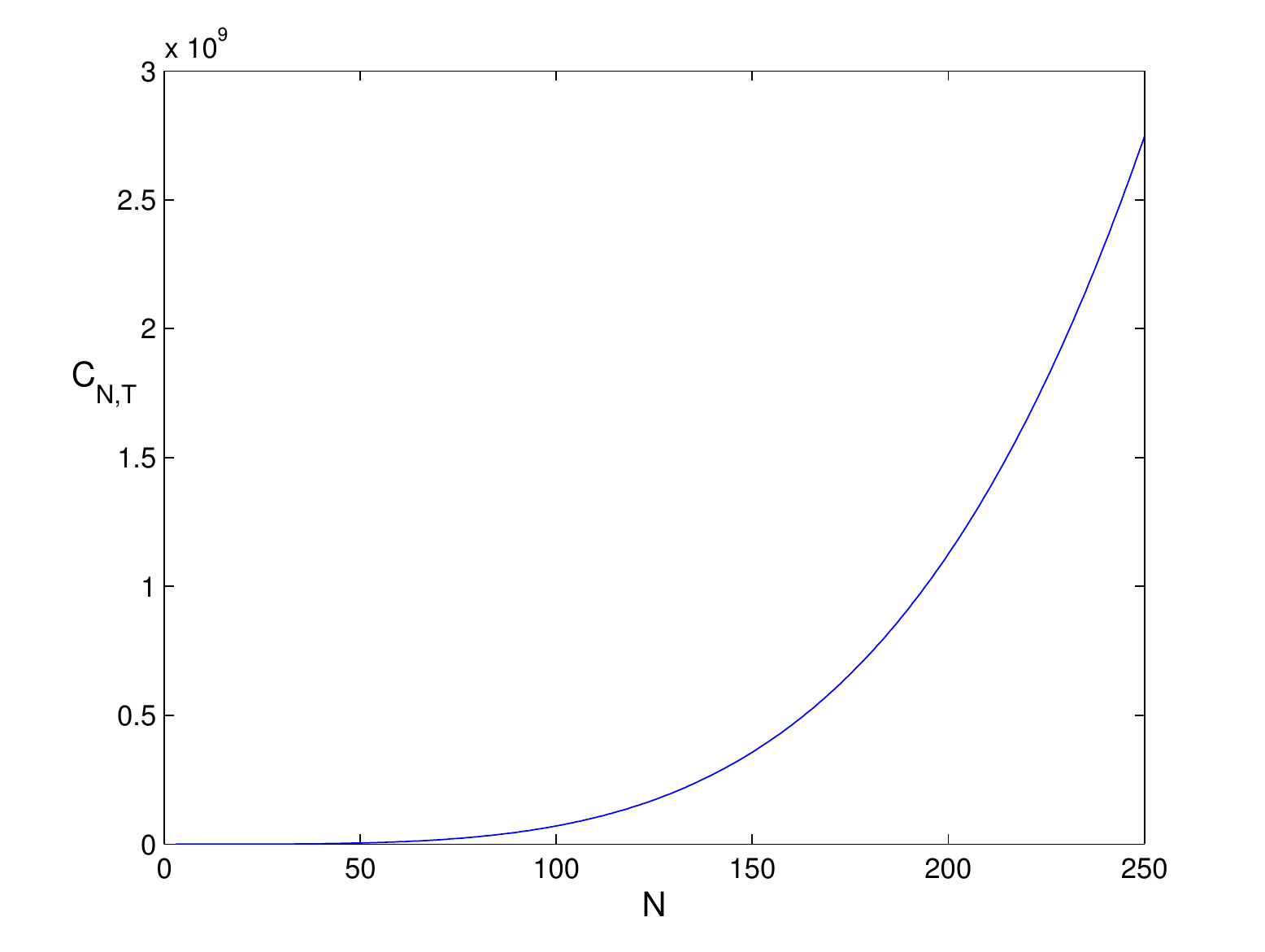}} 
\end{tabular}
\end{center}
\caption{\footnotesize{Values of $c_{N,T}$ (left) and $C_{N,T}$ (right) for the non-symmetric Nitsche's method, without the term $\gamma N^2 u^N(1,t)$  in their definitions (\ref{NitscheNonSymobs}) and (\ref{NitscheNonSymtrace}), with $\gamma=1$ and $T=8$}}\label{comparaisonNitscheNSnonunif}
\end{figure}

According to our numerical investigations, the term $\gamma N^2 u^N(1,t)$ is necessary to obtain a uniform and strictly positive lower bound on $c_{N,T}$ and that removing this term in the observation does not allow one to recover an uniform upper bound on $C_{N,T}$. Figure \ref{comparaisonNitscheNSnonunif} illustrates this in the same conditions ($\gamma=1$, $T=8$) than in the previous example.

\section{A numerical example}

In this section we show the influence that all the remedies (spectral filter, mixed formulation and Nitsche's method) to the non-uniform observability propertiy have on the computation of the associated control problem. 

In this example, we take  

\begin{eqnarray}\label{donnees2}
      y_0(x)&=&x+1 \nonumber \\
      y_1(x)&=&0
\end{eqnarray}  

The corresponding values of the exact control $v$ and the initial data $(u_0,u_1)$ for the adjoint homogeneous wave equation giving the control $v=u_x(1,t)$ are 
\begin{eqnarray*}\label{solutions2}
      u_0(x)&=&0  \\
      u_1(x)&=&-x/4 -1/4  \\
      v(t)&=&-t/4+1/2 \; \mbox{if} \; 0<t<4, \\
      v(t)&=&-t/4+3/2 \; \mbox{if} \, 4<t<8
\end{eqnarray*}  

%
%
%

\begin{figure}[ht]
\begin{center}
\begin{tabular}{c}
 \mbox{\includegraphics[height=4cm]{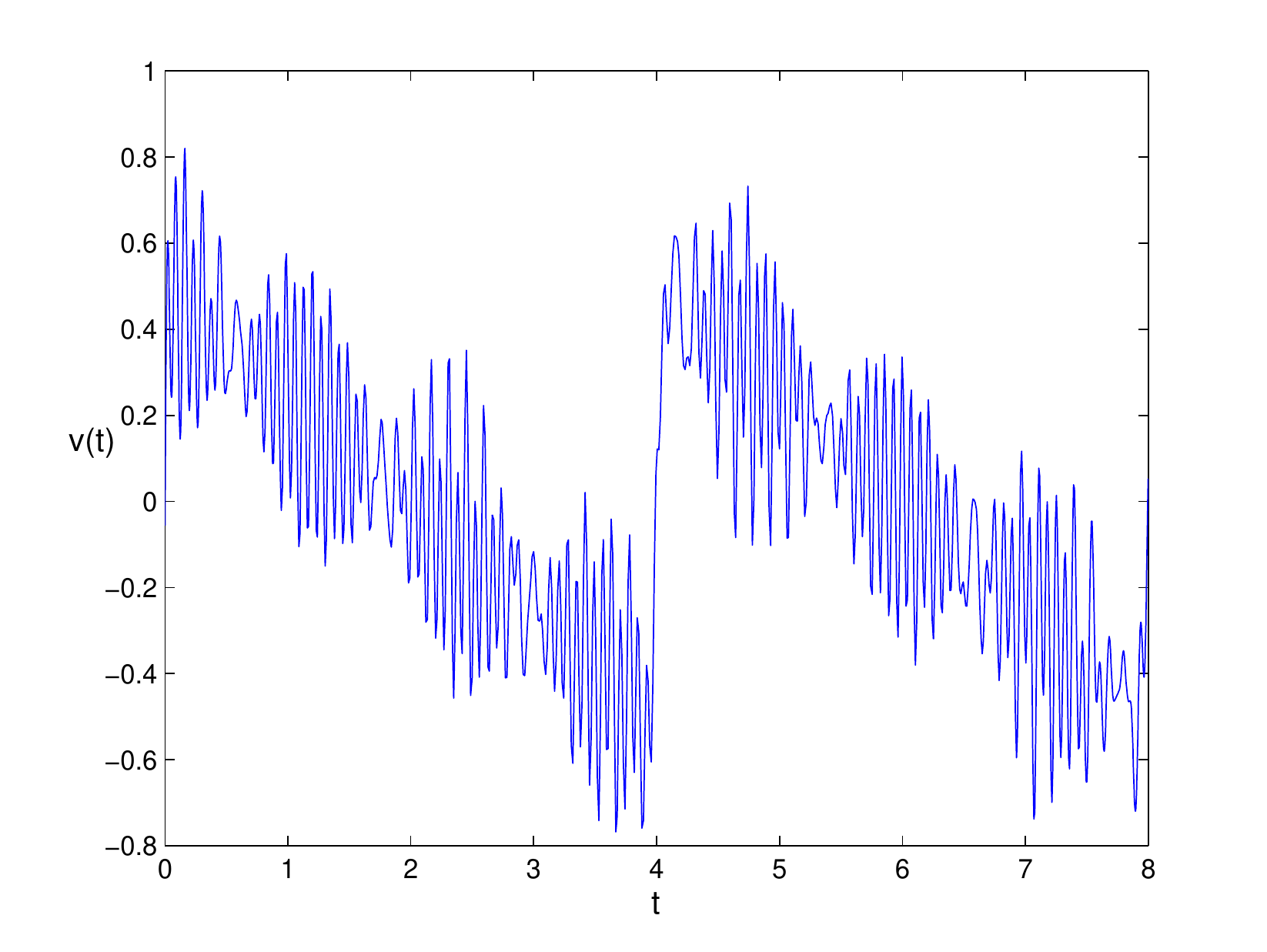}}
 \mbox{\includegraphics[height=4cm]{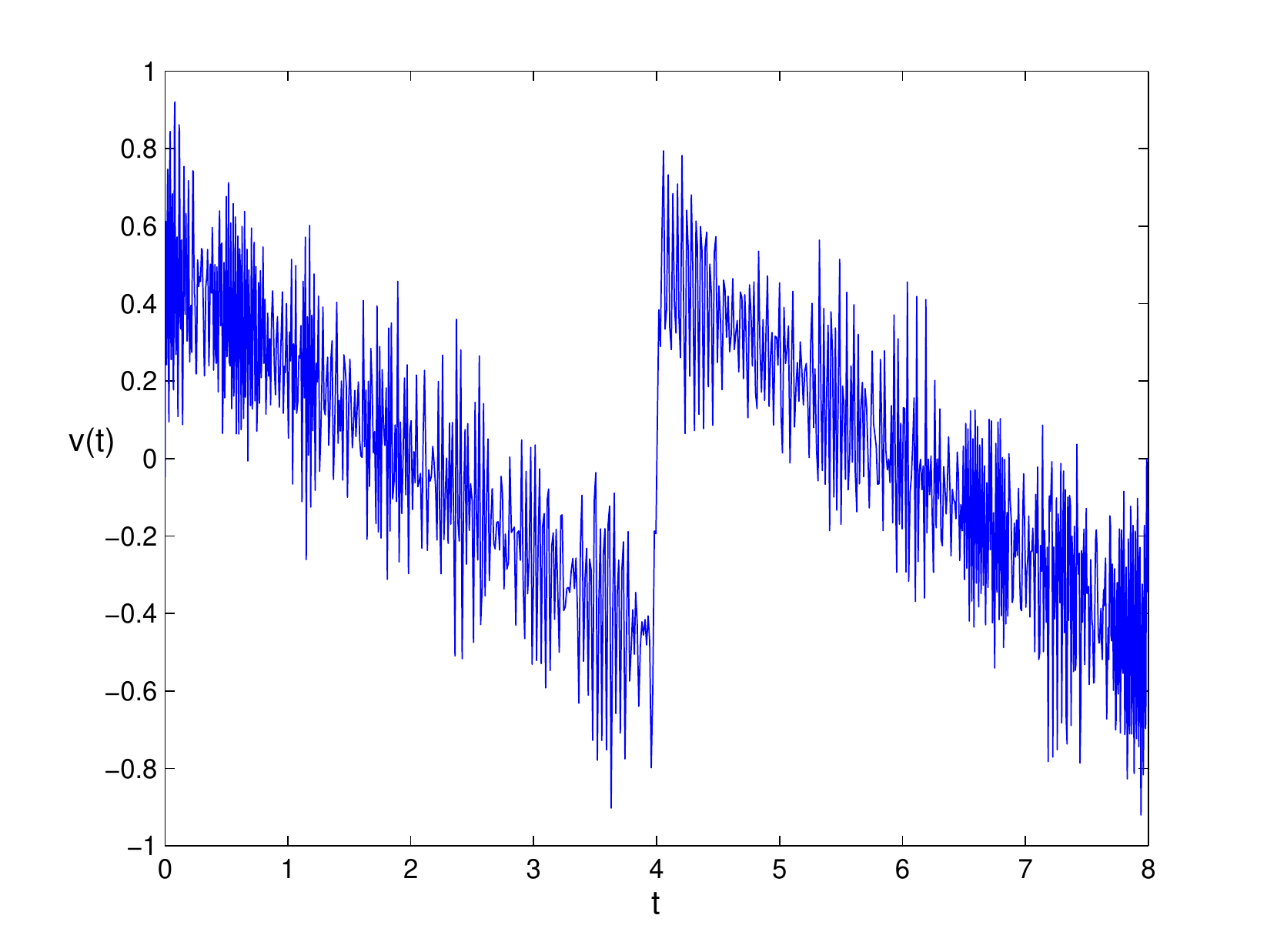}} \\
 \mbox{\includegraphics[height=4cm]{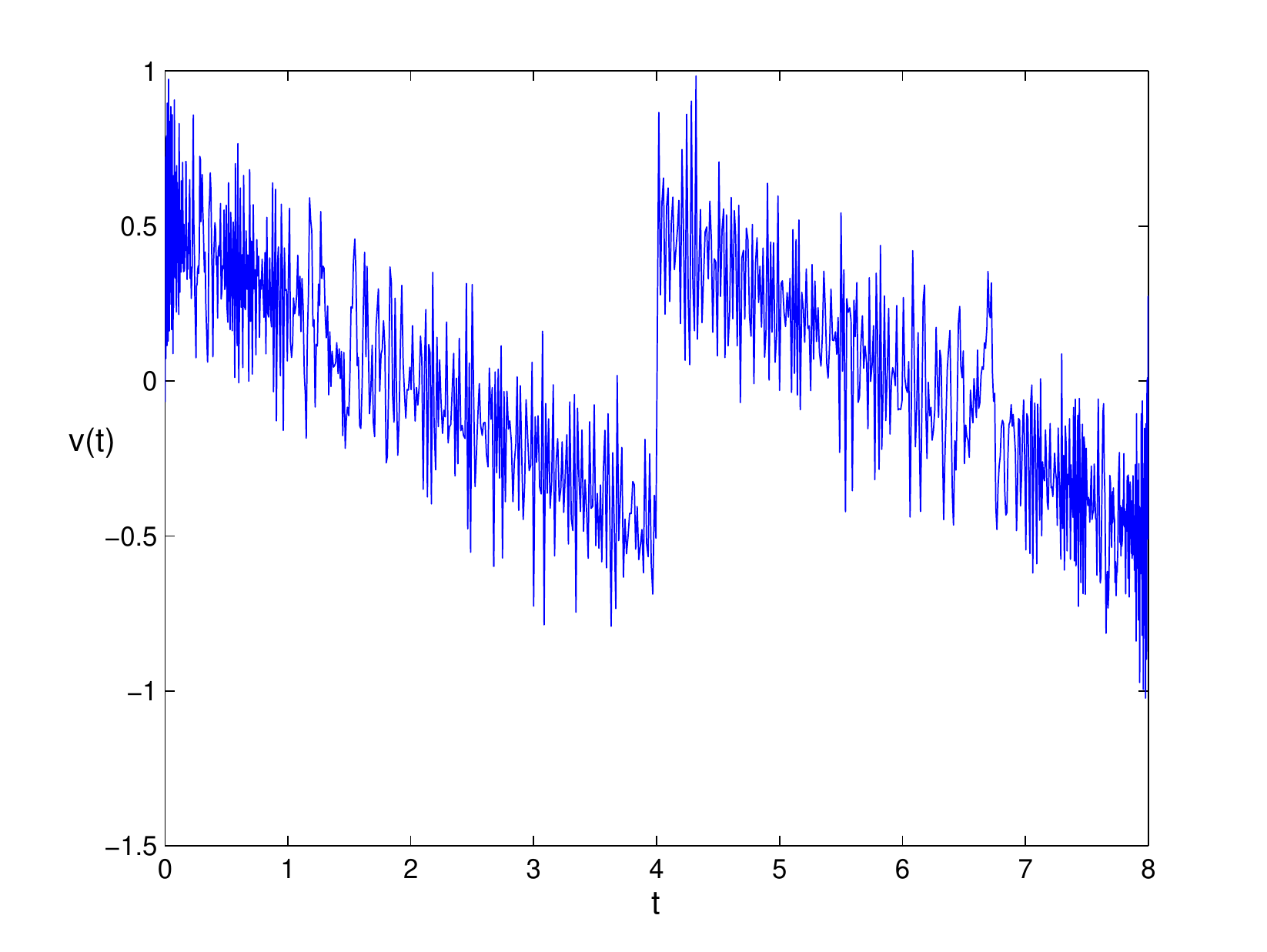}} 
 \mbox{\includegraphics[height=4cm]{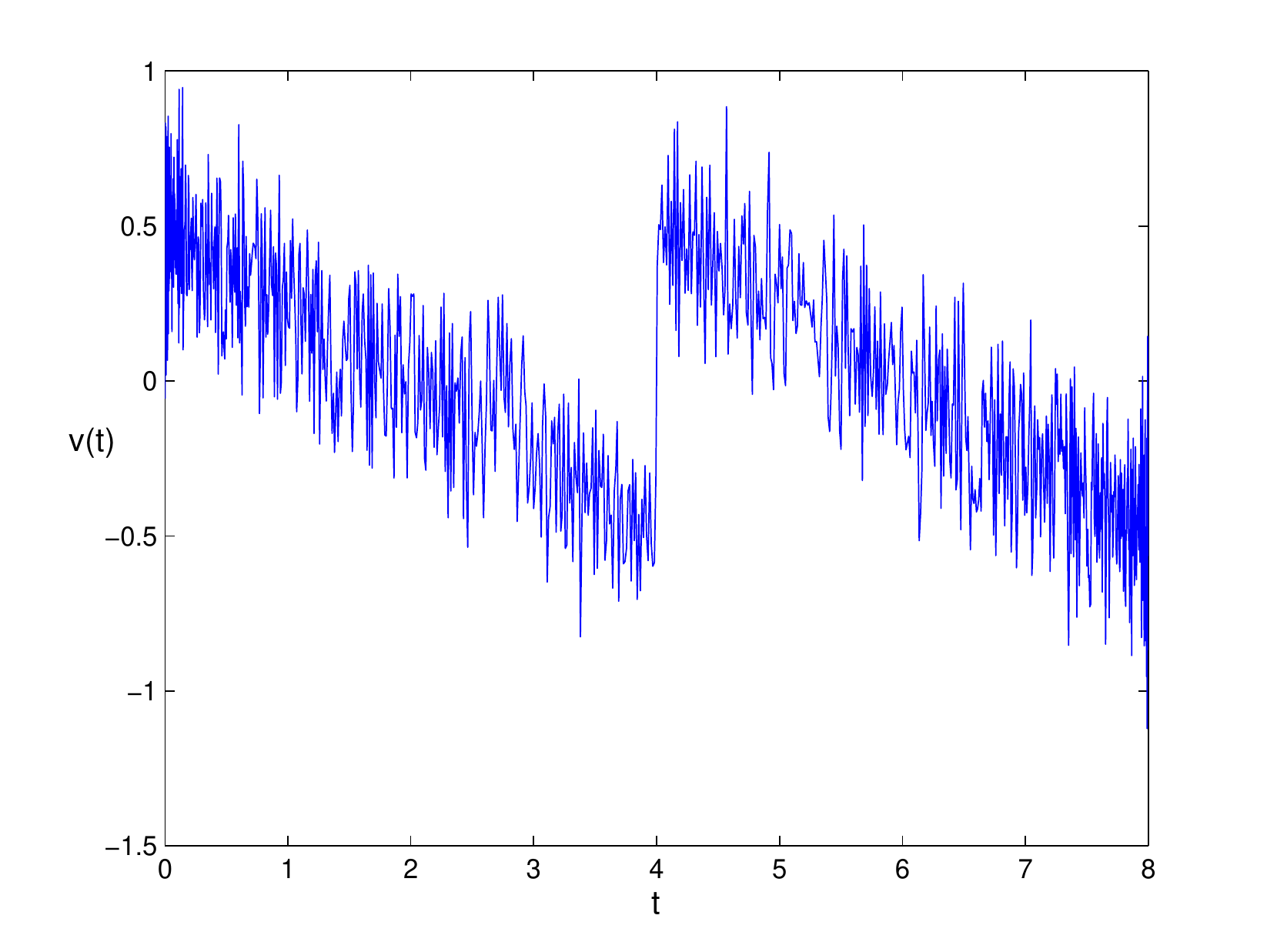}} \\
\end{tabular}
\end{center}
\caption[Contr\^ole num\'erique issu de la formulation (\ref{eqdiscret}) associ\'e aux donn\'ees initiales (\ref{donnees2})]{\footnotesize{Numerical control $v^N(t)$ obtained with the  Legendre Galerkin method (\ref{obsdiscret}) for $N=32$ (top left) $N=64$ (top right), $N=128$ (bottom left) and $N=256$ (bottom right). $T=8$.}}\label{controlenormal}
\end{figure}

In Figure \ref{controlenormal} we show the approximation $v^N(t)$ obtained with the Legendre Galerkin method and without any of the remedies studied here. This approximation is highly oscillatory, even for low values of $N$, and this pathology amplifies with growing values of $N$.    

\begin{figure}[ht]
\begin{center}
\begin{tabular}{c}
 \mbox{\includegraphics[height=4cm]{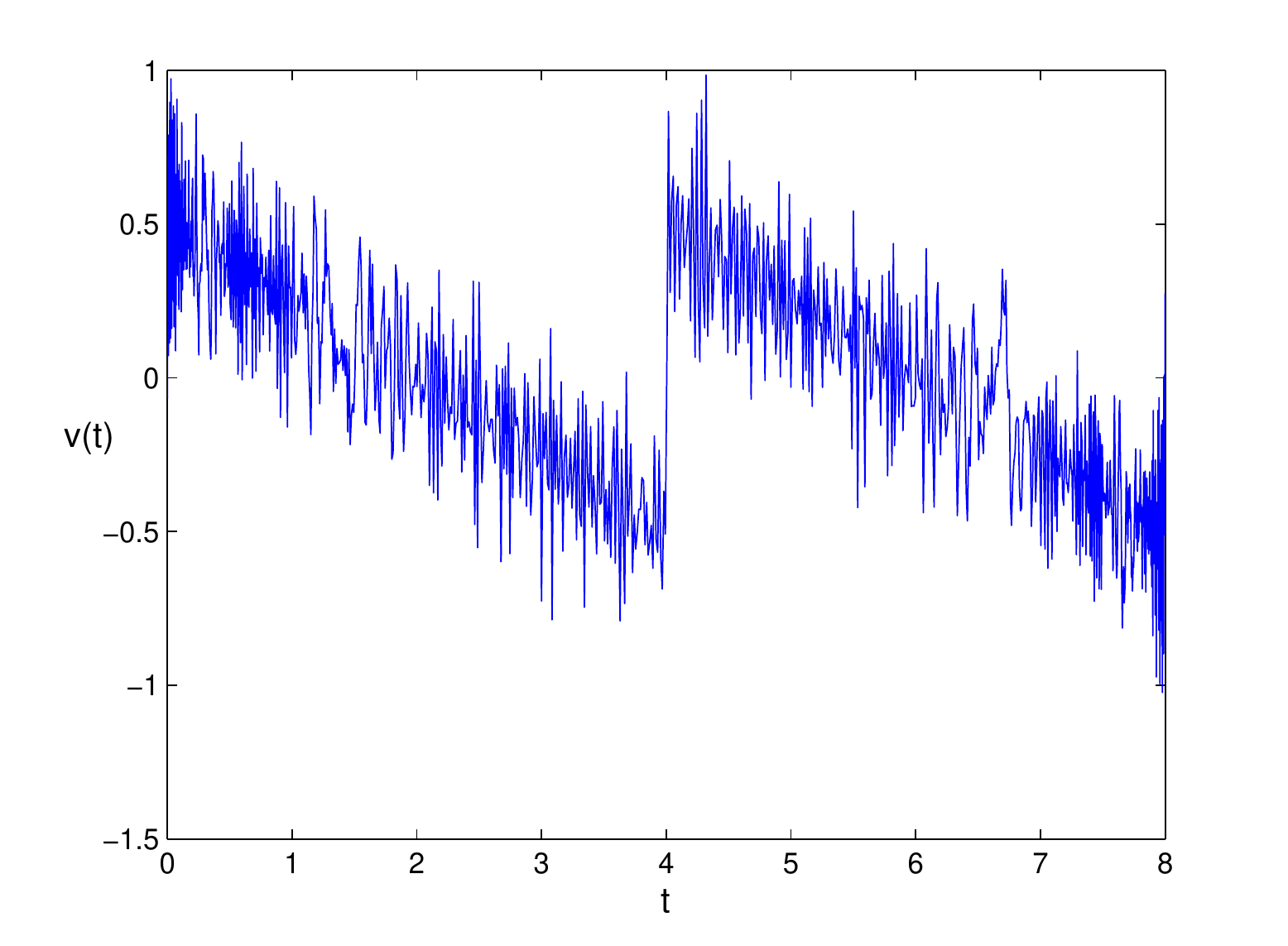}} 
 \mbox{\includegraphics[height=4cm]{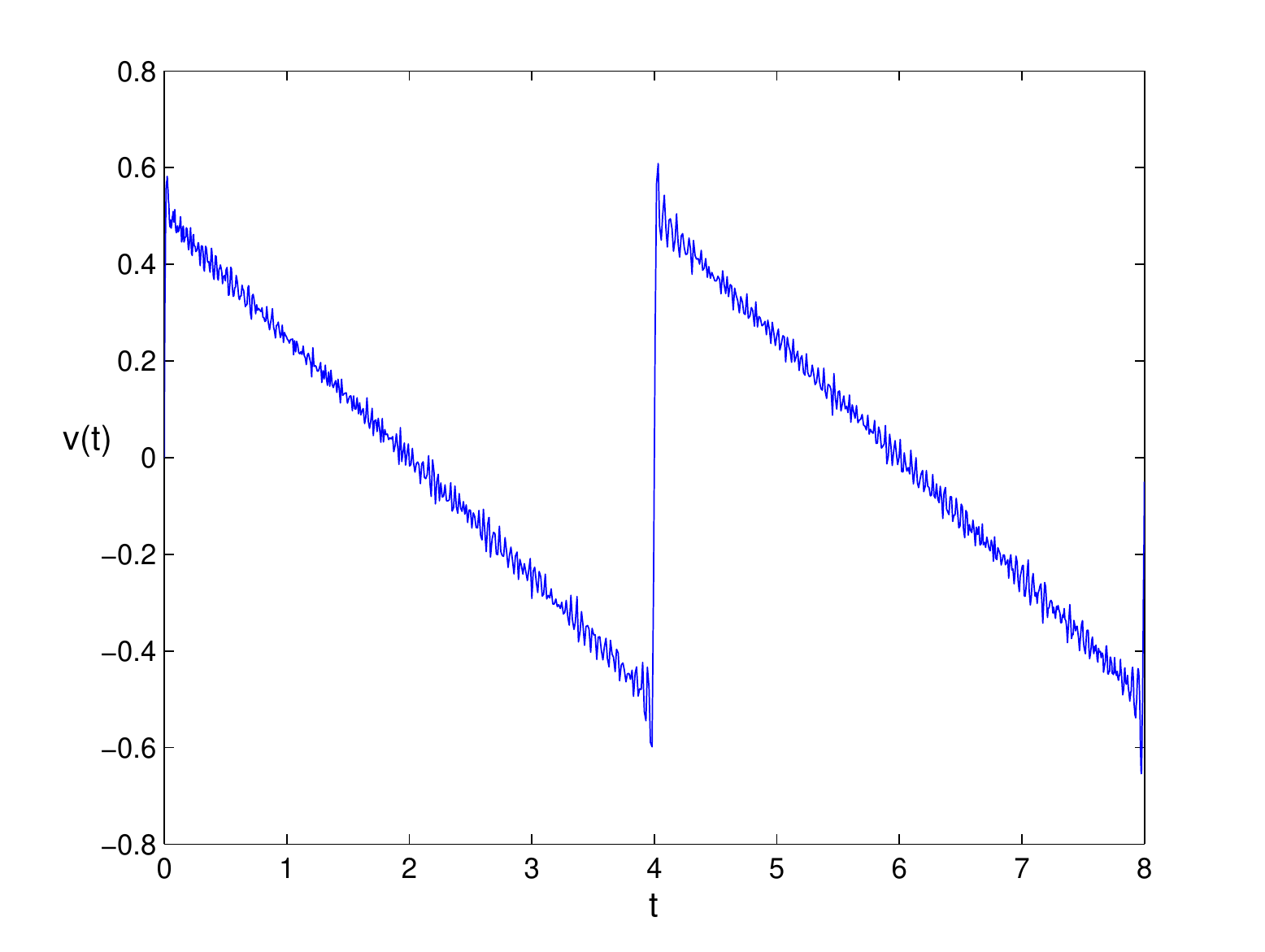}} \\
 \mbox{\includegraphics[height=4cm]{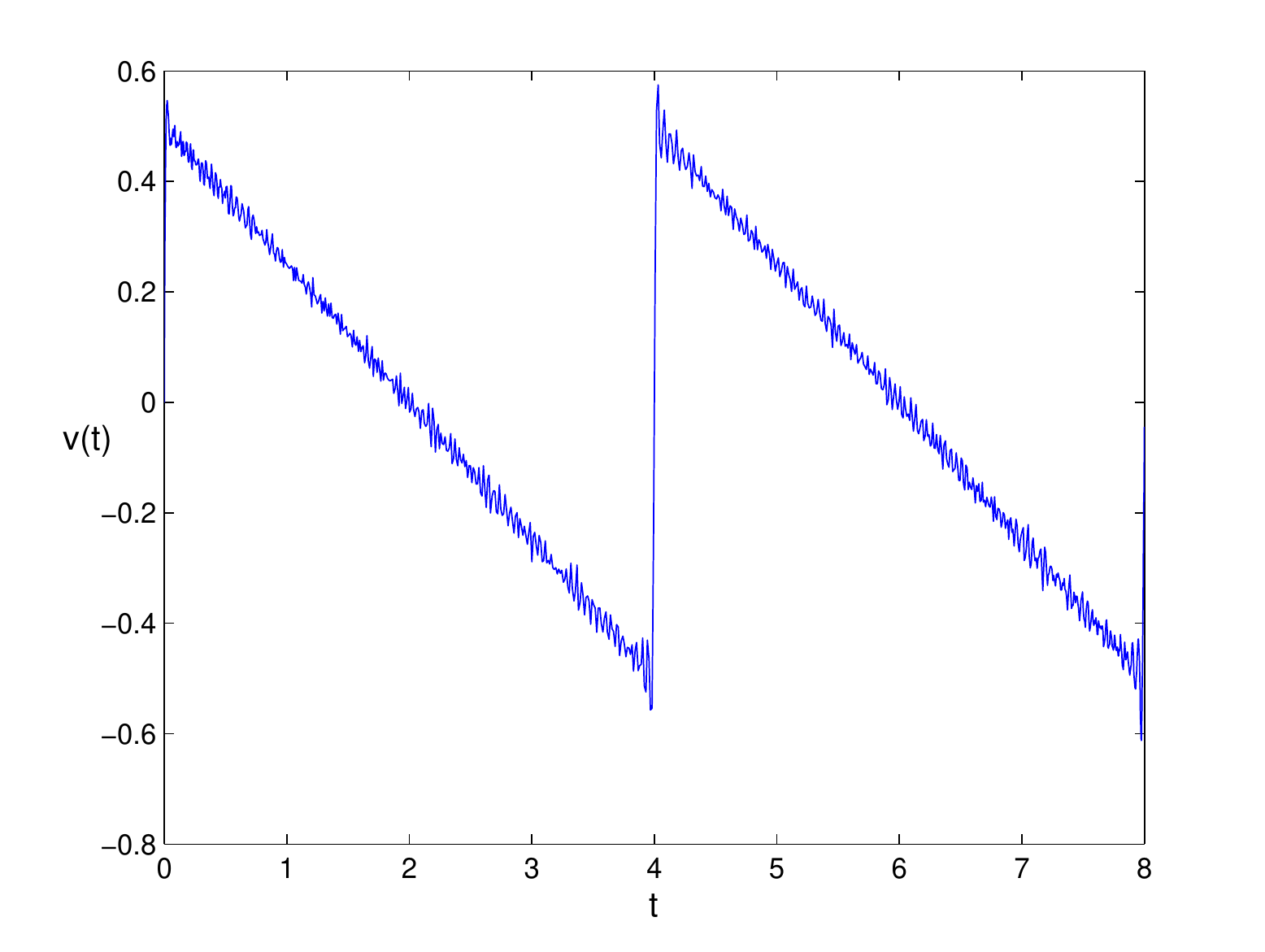}} 
 \mbox{\includegraphics[height=4cm]{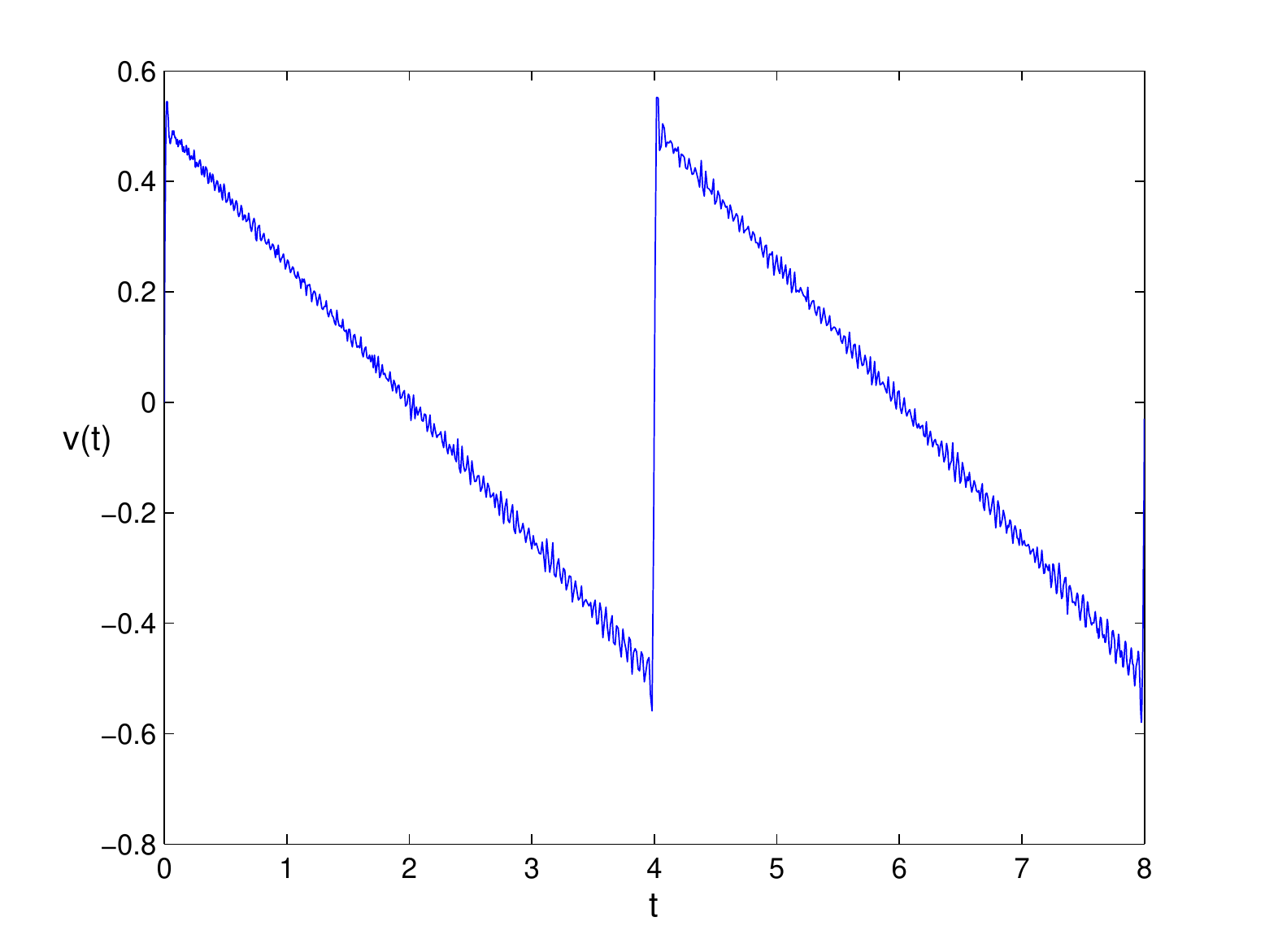}} \\
 \mbox{\includegraphics[height=4cm]{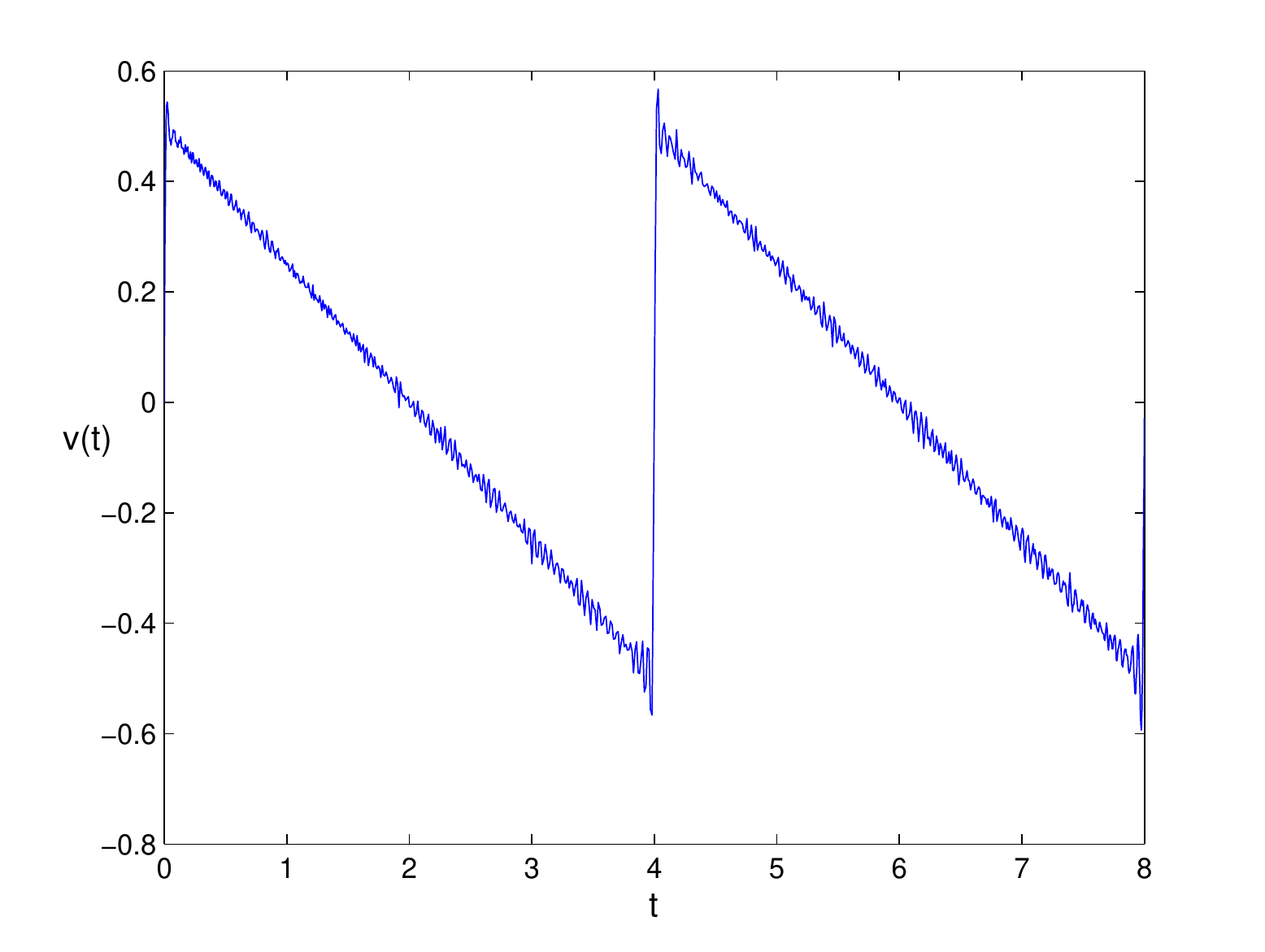}} 
 \mbox{\includegraphics[height=4cm]{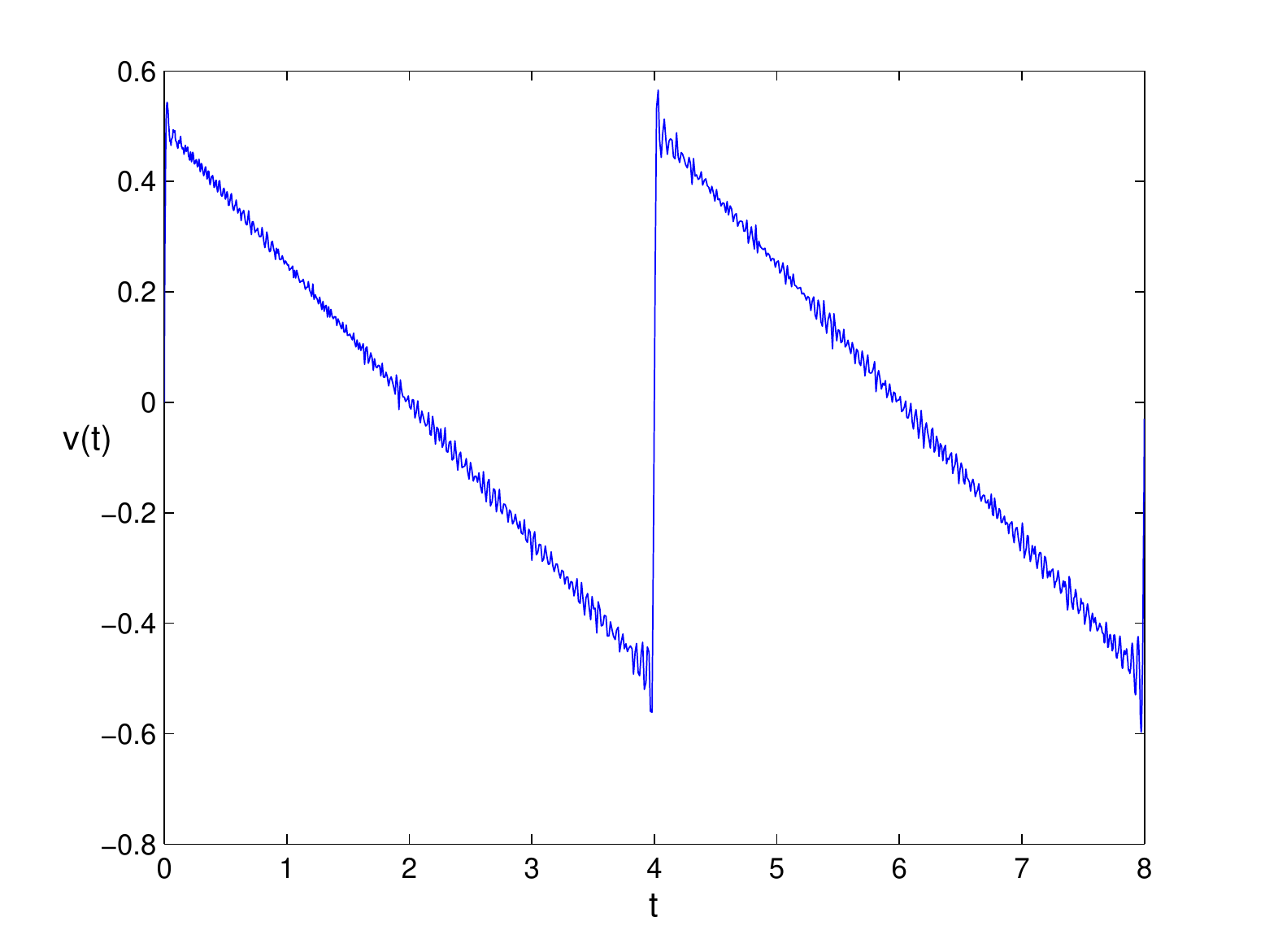}} \\
\end{tabular}
\end{center}
\caption[Comparaison du contr\^ole num\'erique $v^N(t)$ associ\'e \`a (\ref{donnees2}) pour quelques formulations pour $N$ fix\'e \`a 128]{\footnotesize{Numerical controls with $N=128$ with (from top to bottom): the classical formulation (\ref{obsdiscret}) (left), Ces\`aro smoothing (right) - Exponential smoothing (p=6) (left), the mixed formulation (right) - symmetric Nitsche's method  ($\gamma=1$) (left), non-symmetric Nitsche's method  ($\gamma=1$) (right)}}\label{comparaisoncontroles}
\end{figure}

In Figure \ref{comparaisoncontroles} we show the approximations $v^N(t)$, for $N=128$, obtained with the different remedies studied here: spectral filtering with Cesaro and exponential filters, mixed formulation of the equations and Nitsche's method in both its symmetrical and non-symmetrical forms. All these remedies to the non-uniform observability property of the classical Legendre Galerkin method seem to give good approximations of the exact control. There is some  oscillatory pollutions for each of them, but always with small amplitudes. 

\begin{figure}[ht!]
\begin{center}
\begin{tabular}{c}
 \mbox{\includegraphics[height=4cm]{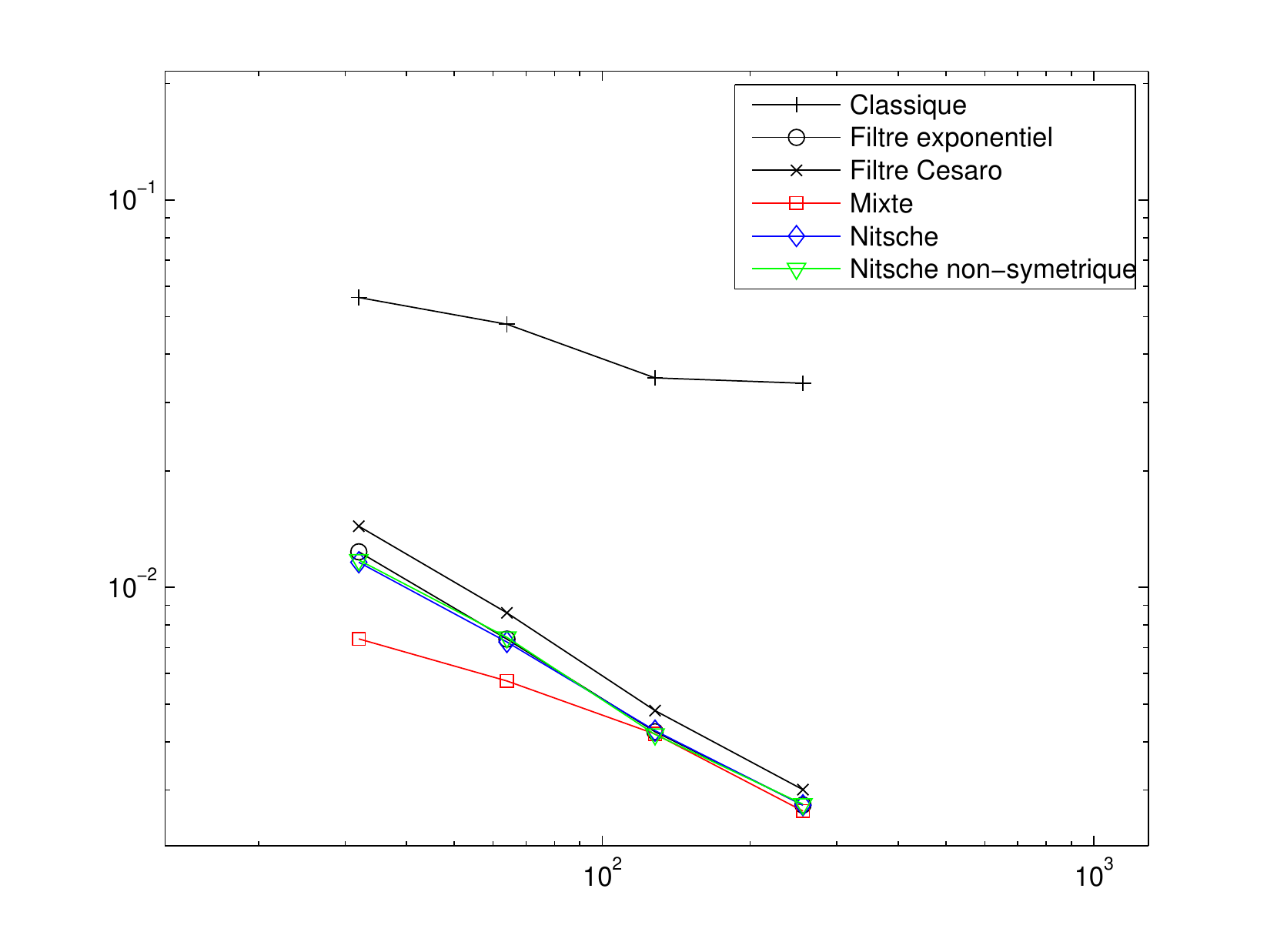}}
 \mbox{\includegraphics[height=4cm]{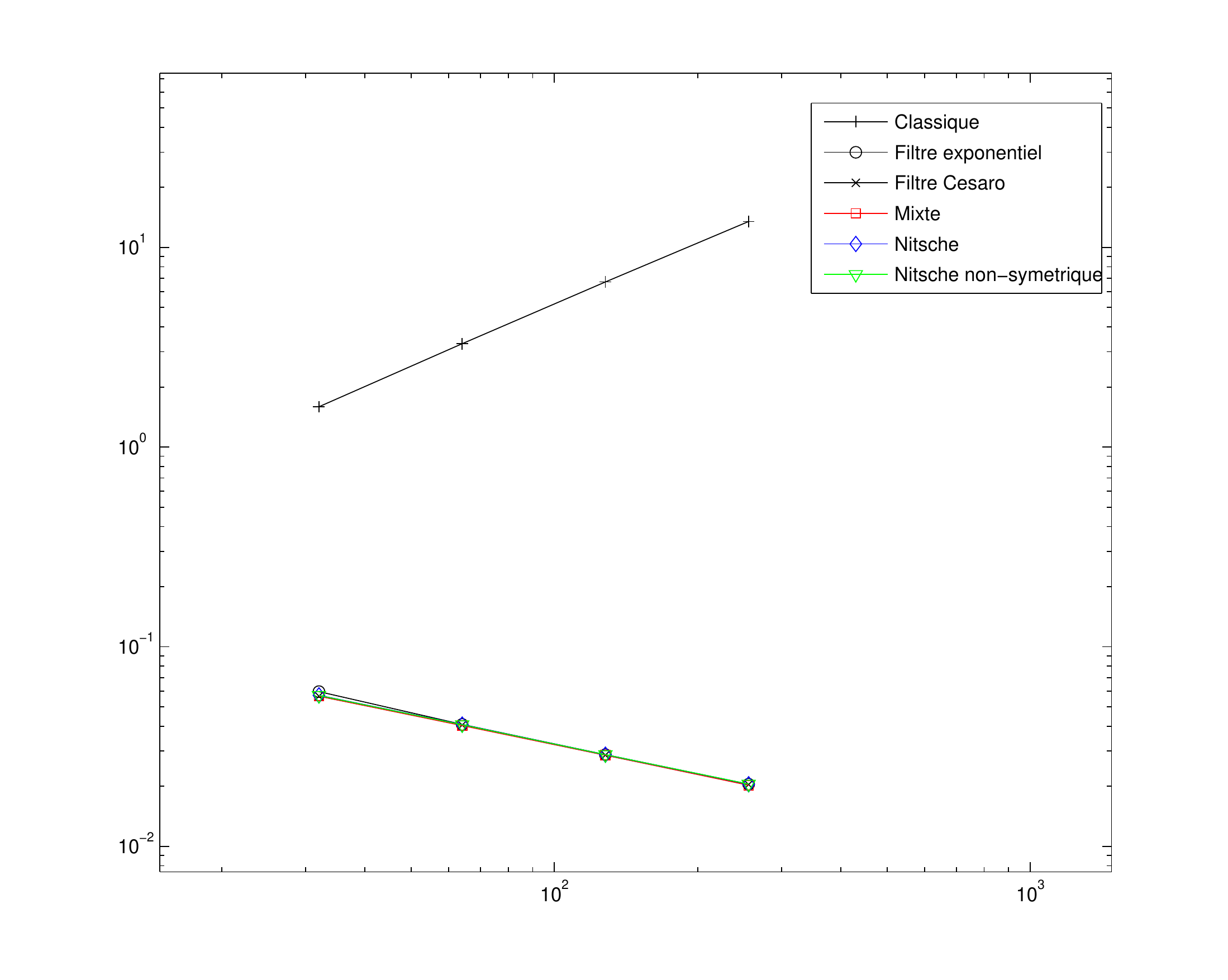}}\\
 \mbox{\includegraphics[height=4cm]{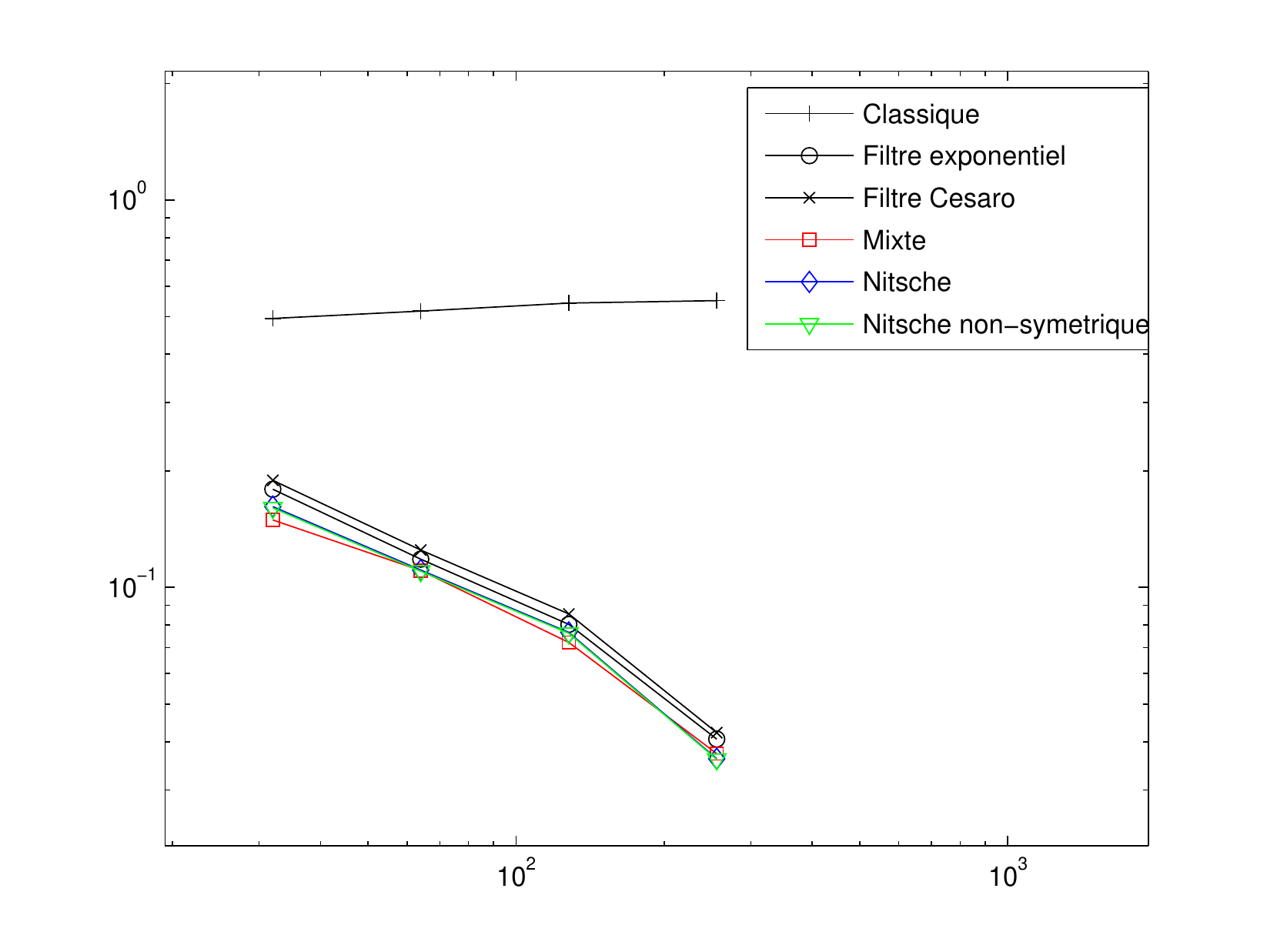}}
\end{tabular}
\end{center}
\caption[Convergence des donn\'ees initiales $u^N_0,u^N_1$ et du contr\^ole $v^N(t)$ associ\'ees \`a (\ref{donnees2}) pour quelques formulations]{\footnotesize{Errors $|u_0^N-u_0|_{H^1_0}$, $\|u_1^N-u_1\|_{L^2}$ and $\|v^N-v\|_{L^2(0,T)}$ for the classical Legendre Galerkin method, and the remedies studied here :  Ces\`aro and exponential (p=6) filters, the mixed formulation and Nitsche's method ($\gamma=1$) in its symmetric and non-symmetric forms for $N=32,\,64,\,128,\,256$}}\label{convergences}
\end{figure}

In order to confirm these good results, in Figure \ref{convergences} we show the behaviour of the errors $|u_0^N-u_0|_{H^1_0}$, $\|u_1^N-u_1\|_{L^2}$ and $\|v^N-v\|_{L^2(0,T)}$ with respect to $N$. For the classical Legendre Galerkin method, without any remedy, there is no convergence of the approximations. However, for each of the remedies studied here, convergence appears to hold, with approximately the same rate. 

\section{Appendix : Computation of $c_{N,T}$ and $C_{N,T}$}

In order to explain how the values of $c_{N,T}$ and $C_{N,T}$ are computed, let's consider the first order matrix form of a semi-discretization of (\ref{obs}),
\begin{eqnarray*}
\left( \begin{array}{c}
        \mathbf{a}_N'(t) \\
        \mathbf{a}_N''(t)
       \end{array}
   \right)&=&\left( \begin{array}{cc}
        0 & I_{N-1\times N-1} \\
        -M_N^{-1}K_N & 0
       \end{array}
   \right)   \left( \begin{array}{c}
        \mathbf{a}_N(t) \\
        \mathbf{a}_N'(t)
       \end{array}
   \right)   \\
   &=:&A_N
   \left( \begin{array}{c}
        \mathbf{a}_N(t) \\
        \mathbf{a}_N'(t)
       \end{array}
   \right).
\end{eqnarray*}
The solution is 
\begin{equation*}
\left( \begin{array}{c}
        \mathbf{a}_N(t) \\
        \mathbf{a}_N'(t)
       \end{array}
   \right) =e^{A_Nt}
   \left( \begin{array}{c}
        \mathbf{a}_N(0) \\
        \mathbf{a}_N'(0)
       \end{array}
   \right).
\end{equation*}  
Let's also write the observation in matrix form :
\begin{equation*}
u^N_x(1,t)=\left(C_{N}^t,\left( \begin{array}{c} \mathbf{a}_N(t) \\ \mathbf{a}_N'(t) \end{array} \right)\right)_{\RR^{2N}}.
\end{equation*}
We then have
\begin{eqnarray*}
\int_0^T\left|u^N_x(1,t)\right|^2 \, \textrm{dt} &=& \left(\left( \int_0^T e^{A_N^tt}C_N^tC_Ne^{At}\, \textrm{dt} \right)\left( \begin{array}{c} \mathbf{a}_N(0) \\ \mathbf{a}_N'(0) \end{array} \right), \left( \begin{array}{c} \mathbf{a}_N(0) \\ \mathbf{a}_N'(0) \end{array} \right)\right)_{\RR^{2N}}\\
&=:& \left(W^N(0,T) \left( \begin{array}{c} \mathbf{a}_N(0) \\ \mathbf{a}_N'(0) \end{array} \right), \left( \begin{array}{c} \mathbf{a}_N(0) \\ \mathbf{a}_N'(0) \end{array} \right)\right)_{\RR^{2N}}.
\end{eqnarray*}
On the other hand, the discrete energy writes
\begin{equation*}
E(u^N(0))=\dfrac{1}{2}\left(\left( \begin{array}{cc}
                                        K_N & 0 \\ 
                                        0 & M_N
                                       \end{array}\right)\left( \begin{array}{c} \mathbf{a}_N(0) \\ \mathbf{a}_N'(0) \end{array} \right),\left( \begin{array}{c} \mathbf{a}_N(0) \\ \mathbf{a}_N'(0) \end{array} \right)\right)_{\RR^{2N}}.
\end{equation*}
The values of $c_{N,T}$ and $C_{N,T}$ are then, respectively, the minimal and the maximal eigenvalues of the generalized eigenproblem  : find $(\lambda,U)\in \RR\times\RR^{2N}$ such that 
\begin{equation}\label{eigenvaluesAnnexe}
W^N(0,T)\,U=\dfrac{\lambda}{2} \left( \begin{array}{cc} 
                                        K_N & 0 \\ 
                                        0 & M_N
                                       \end{array}
 \right)\,U
\end{equation}
The following lemma (\cite{Chen}) gives a way to compute $W^N(0,T)$. 
\begin{lem}
Let 
\begin{eqnarray*}
F(t)&=&\left( \begin{array}{cc}
		  F_{11}(t) & F_{12}(t) \\ 
		  0 & F_{22}(t)
		  \end{array}\right)\\
    &=&\exp\left[ \left( \begin{array}{cc}
		  A_{11} & A_{12} \\ 
		  0 & A_{22}
		  \end{array}\right)t \right].		  
\end{eqnarray*}
Then, 
\begin{equation*}
 F_{11}(t)=e^{A_{11}t}, F_{22}(t)=e^{A_{22}t} \textrm{ et } F_{12}(t)=\int_0^t e^{A_{11}(t-s)}A_{12}e^{A_{22}s} \, \textrm{ds}.
\end{equation*}
\end{lem}
Choosing $A_{11}=-A_N^t$, $A_{12}=C_N^tC_N$ and $A_{22}=A_N$, we have that  
\begin{equation*}
W(0,T)=F_{22}(T)^tF_{12}(T). 
\end{equation*}


\bibliographystyle{plain}
\bibliography{biblio}

\begin{thebibliography}{10}

\bibitem{ARN}
D.~N. Arnold.
\newblock An interior penalty finite element method with discontinuous
  elements.
\newblock {\em SIAM J. Num. Anal.}, 19(4):742--760, 1982.

\bibitem{JU}
T.~Z. Boulmezaoud and J.~M. Urquiza.
\newblock On the eigenvalues of the spectral second order differentiation
  operator and application to the boundary observability of the wave equation.
\newblock {\em J. Sci. Comput.}, 31(3):307--345, 2007.

\bibitem{CAN1}
C.~Canuto, M.~Y. Hussaini, A.~Quarteroni, and T.~A. Zang.
\newblock {\em Spectral methods}.
\newblock Scientific Computation. Springer-Verlag, Berlin, 2006.
\newblock Fundamentals in single domains.

\bibitem{MICU}
C.~Castro and S.~Micu.
\newblock Boundary controllability of a linear semi-discrete 1-{D} wave
  equation derived from a mixed finite element method.
\newblock {\em Numer. Math.}, 102(3):413--462, 2006.

\bibitem{Chen}
T.~Chen and B.~Francis.
\newblock {\em Optimal sampled-data control systems}.
\newblock Communications and Control Engineering Series. Springer-Verlag London
  Ltd., London, 1996.
\newblock Reprint of the 1995 original.

\bibitem{ErvedozaZuazua2012}
S.~Ervedoza and E.~Zuazua.
\newblock {\em The wave equation : control and numerics}.
\newblock in Control of Partial Differential Equations, P.M. Cannarsa and J.M.
  Coron Eds., 245-340. Springer, New York, 2012.

\bibitem{ErvedozaZuazua2013}
S.~Ervedoza and E.~Zuazua.
\newblock {\em Numerical approximation of exact controls for waves}.
\newblock Springer Briefs in Mathematics. Springer, New York, 2013.

\bibitem{GLO}
R.~Glowinski.
\newblock Ensuring well-posedness by analogy: {S}tokes problem and boundary
  control for the wave equation.
\newblock {\em J. Comput. Phys.}, 103(2):189--221, 1992.

\bibitem{GLL}
R.~Glowinski, C.~H. Li, and J.-L. Lions.
\newblock A numerical approach to the exact boundary controllability of the
  wave equation. {I}. {D}irichlet controls: description of the numerical
  methods.
\newblock {\em Japan J. Appl. Math.}, 7(1):1--76, 1990.

\bibitem{GROTE}
M.~J. Grote, A.~Schneebeli, and D.~Sch{\"o}tzau.
\newblock Discontinuous {G}alerkin finite element method for the wave equation.
\newblock {\em SIAM J. Numer. Anal.}, 44(6):2408--2431 (electronic), 2006.

\bibitem{hansbo}
P.~Hansbo.
\newblock Nitsche's method for interface problems in computational mechanics.
\newblock {\em GAMM-Mitt.}, 28(2):183--206, 2005.

\bibitem{HEST}
J.~S. Hesthaven and R.~M. Kirby.
\newblock Filtering in {L}egendre spectral methods.
\newblock {\em Math. Comp.}, 77(263):1425--1452, 2008.

\bibitem{infantezuazua}
J.A. Infante and E.~Zuazua.
\newblock Boundary observability for the space semi-discretization of the 1-d
  wave equation.
\newblock {\em Math. Model. Num. Anal.}, 33:407--438, 1999.

\bibitem{KOMOR}
V.~Komornik.
\newblock {\em Exact controllability and stabilization}.
\newblock RAM: Research in Applied Mathematics. Masson, Paris, 1994.
\newblock The multiplier method.

\bibitem{JL88}
J.-L. Lions.
\newblock {\em Contr\^olabilit\'e exacte, perturbations et stabilisation de
  syst\`emes distribu\'es. {T}ome 1}, volume~8 of {\em Recherches en
  Math\'ematiques Appliqu\'ees [Research in Applied Mathematics]}.
\newblock Masson, Paris, 1988.
\newblock Contr{\^o}labilit{\'e} exacte. [Exact controllability], With
  appendices by E. Zuazua, C. Bardos, G. Lebeau and J. Rauch.

\bibitem{MaricaZuazua}
A.~Marica and E.~Zuazua.
\newblock {\em Symmetric Discontinuous Galerkin Methods for 1-D waves. Fourier
  Analysis, Propagation, Observability, and Applications.}
\newblock Springer Briefs in Mathematics. Springer, New York, 2014.

\bibitem{NEGZ}
M.~Negreanu and E.~Zuazua.
\newblock Convergence of a multigrid method for the controllability of a 1-d
  wave equation.
\newblock {\em C. R. Math. Acad. Sci. Paris}, 338(5):413--418, 2004.

\bibitem{Nitsche}
J.~Nitsche.
\newblock \"{U}ber ein {V}ariationsprinzip zur {L}\"osung von
  {D}irichlet-{P}roblemen bei {V}erwendung von {T}eilr\"aumen, die keinen
  {R}andbedingungen unterworfen sind.
\newblock {\em Abh. Math. Sem. Univ. Hamburg}, 36:9--15, 1971.
\newblock Collection of articles dedicated to Lothar Collatz on his sixtieth
  birthday.

\bibitem{SHEN}
J.~Shen.
\newblock Efficient spectral-{G}alerkin method. {I}. {D}irect solvers of
  second- and fourth-order equations using {L}egendre polynomials.
\newblock {\em SIAM J. Sci. Comput.}, 15(6):1489--1505, 1994.

\bibitem{sontag}
E.~Sontag.
\newblock {\em Mathematical Control Theory}.
\newblock Springer, New York, 1998.

\bibitem{VANDEVEN}
H.~Vandeven.
\newblock On the eigenvalues of second-order spectral differentiation
  operators.
\newblock {\em Comput. Methods Appl. Mech. Engrg.}, 80(1-3):313--318, 1990.

\bibitem{WAR}
T.~Warburton and J.~S. Hesthaven.
\newblock On the constants in {$hp$}-finite element trace inverse inequalities.
\newblock {\em Comput. Methods Appl. Mech. Engrg.}, 192(25):2765--2773, 2003.

\bibitem{ZuazuaJMAP}
E.~Zuazua.
\newblock Boundary observability for the finite-difference space
  semi-discretizations of the 2-d wave equation in the square.
\newblock {\em J. Math. Pures Appl.}, 78(5):523--563, 1999.

\bibitem{ZUA1}
E.~Zuazua.
\newblock Propagation, observation, and control of waves approximated by finite
  difference methods.
\newblock {\em SIAM Review}, 47(2):197--243, 2005.

\end{thebibliography}

\end{document}